\documentclass[11pt]{article}
\usepackage[T1]{fontenc} 
\usepackage[utf8]{inputenc}

\usepackage{lmodern}

\usepackage{amssymb} 
\usepackage{latexsym}
\usepackage{epsfig}
\usepackage{lscape,graphicx,color,psfrag}
\usepackage{calrsfs}
\usepackage[cmex10]{amsmath}
\usepackage{amsthm}
\usepackage{natbib}

\theoremstyle{plain}

\newtheorem{assumption}{Assumption}

\newtheorem{theorem}{Theorem}

\theoremstyle{definition}
\newtheorem{definition}{Definition}
\theoremstyle{remark}

\def\argmax{\mathop{\rm arg\,max}}

\newcommand{\Prob}{\ensuremath{\mathsf{P}}}
\newcommand{\Ex}{\ensuremath{\mathsf{E}}}
\newcommand{\conv}{\ensuremath{\mathrm{conv}}} 

\usepackage{appendix}

\usepackage[top=1in, bottom=1in, left=0.9in, right=0.9in]{geometry}

\usepackage[colorlinks=true,urlcolor=blue,citecolor=blue,linkcolor=blue]{hyperref}
\usepackage[noabbrev,nameinlink]{cleveref}     

\def\EMAIL#1{\href{mailto:#1}{#1}}

\begin{document}

\title{Dynamic priority allocation via restless bandit marginal
  productivity indices}


\author{Jos\'e Ni\~no-Mora 
\\ Department of Statistics \\
    Carlos III University of Madrid  \\  \EMAIL{jose.nino@uc3m.es}, \href{http://alum.mit.edu/www/jnimora}{http://alum.mit.edu/www/jnimora} \\
      ORCID: \href{http://orcid.org/0000-0002-2172-3983}{0000-0002-2172-3983}}
 
\date{Published in \textit{TOP}, vol.\  15, pp.\ 161--198, 2007 \\ \vspace{.1in}
DOI: \href{https://doi.org/10.1007/s11750-007-0025-0}{10.1007/s11750-007-0025-0}}

\maketitle


\begin{abstract}%
This paper surveys recent work by the author on the theoretical and algorithmic aspects of restless bandit indexation as well as on its application to a variety of problems involving the dynamic allocation of priority to multiple stochastic projects. The main aim is to present ideas and methods in an accessible form that can be of use to researchers addressing problems of such a kind. Besides building on the rich literature on bandit problems, our approach draws on ideas from linear programming, economics, and multi-objective optimization. In particular, it was motivated to address issues raised in the seminal work of Whittle (Restless bandits: activity allocation in a changing world. In: Gani J. (ed.) A Celebration of Applied Probability, J. Appl. Probab., vol. 25A, Applied Probability Trust, Sheffield, pp. 287-298, 1988) where he introduced the index for restless bandits that is the starting point of this work. Such an index, along with previously proposed indices and more recent extensions, is shown to be unified through the intuitive concept of ``marginal productivity index'' (MPI), which measures the marginal productivity of work on a project at each of its states. In a multi-project setting, MPI policies are economically sound, as they dynamically allocate higher priority to those projects where work appears to be currently more productive. Besides being tractable and widely applicable, a growing body of computational evidence indicates that such index policies typically achieve a near-optimal performance and substantially outperform benchmark policies derived from conventional approaches.
\end{abstract}%


\tableofcontents
\newpage

\section{Introduction}
\label{s:intro}
The overarching concern with making best use of that most precious
 resource, time,
leads us to ponder how to set priorities among the multifarious
activities vying for our attention.
Thus, we must decide over time whether to engage in projects of
potentially high reward yet unlikely success, or on less
rewarding but more realistic alternatives, revising priorities over time in
light of actual progress and future prospects. 
Similar issues arise in the automatic control of modern technological
systems, such as those in 
manufacturing and computer-communication networks, where the flow
of distinct traffic streams can be regulated
 by dynamically prioritizing access to shared resources such as machines or
 transmission channels.
The high level of discretionarity allowed in such decisions,
as well as their often substantial impact on system performance, raises the
 possibility of  optimizing the latter through appropriate
 design of the priority policy adopted. 

Yet, while many such problems are readily formulated in the framework of 
Markov decision processes (MDPs), their computational solution via the
conventional dynamic programming (DP) technique is typically
intractable, due to the well-known curse of dimensionality. 
As for analytical solutions, they are only available for a few models
under rather special conditions.
Such a state of affairs motivates investigation of heuristic policies
that, while not optimal, are both tractable and come close to achieving desired performance objectives.

Perhaps the most natural and simple class of priority allocation policies is based on
use of \emph{priority indices}. 
Thus, if one must dynamically prioritize work on multiple stochastic projects, an
\emph{index} is defined for each as a function of its state. 
The resultant priority-index policy engages at
each time the required number of projects with currently larger index values.
Yet, such a class of policies is still overwhelmingly large, which
motivates the quest for ideas that guide us to design sound priority
indices yielding good or even optimal policies.

The earliest result on optimality of a priority-index rule is given in
\citet{smith56}, which addresses the problem of sequencing a batch of
jobs having known, deterministic processing times and linear holding
costs. A sound priority index for a job in such a setting is given by
the ratio of its holding cost rate per unit time to its
processing time, which measures the rate of cost reduction per unit of
effort expended.
Hence, such an index can be interpreted as a measure of the
\emph{average productivity of work} on the job.
Smith showed that, in the single-machine case, the total weighted
completion time is minimized by
the resultant index rule. 
The optimality of the \emph{Smith index} rule was extended in
\citet{rothk66} to a
model where job durations are stochastic.
 \citet{coxsmt} showed that such a rule also yields an average-optimal policy
 for scheduling a multiclass single-server queue with linear
holding costs. 
A more complex optimal index rule for the latter model's extension
that incorporates Bernoulli feedback between job classes was obtained
by  \citet{kl}, attaching a constant index to each class.

While such index rules are \emph{static}, in that the index of a job
is constant, in other problems researchers have identified optimal \emph{dynamic}
index rules.
In such a vein, 
the seminal, independent work of \citet{sevcik74} and \citet{gijo74}
 stands out. 
Of particular relevance to this paper is the celebrated result in the
latter paper on the optimal solution of the classic
\emph{multiarmed bandit problem} by a dynamic index rule. 
The problem concerns the sequential allocation of work to a collection of stochastic projects modeled as 
Markov chains that, when engaged, yield rewards and change
state. 
The problem is to decide which project to engage at each time to
maximize the expected total discounted reward earned over an infinite
horizon.
The optimal policy turns out to be the priority-index rule
corresponding to the \emph{Gittins
  index}, which measures the maximum rate of expected discounted
reward per unit of expected discounted time that can be achieved under
stopping rules for each initial project state. 
See \citet{gi79,gi89}. 
Hence, again, the ``right'' index is a measure of the average
productivity of work on a project.
For alternative proofs of such a fundamental
 result offering complementary insights
see, e.g., \citet{whit80}, \citet{vawabu}, \citet{we92}, and
\citet{beni}.

Actually, the roots of such a result can be traced to earlier work.
Thus, the Gittins index extends to a general Markovian
setting the index introduced by
 \citet{bellman56} to solve the special Bayesian Bernoulli one-armed
 bandit problem via calibration.

In turn, Bellman drew on the earlier work of 
\citet{bjk56}, as he acknowledged referring to an unpublished version
of that paper, which we regard as the origin of bandit indexation.
\citet[Sect.\ 4]{bjk56}  addresses the problem of 
optimal sequential design of an experiment where one wishes to
maximize the sum of $n$ observations,  to be chosen
sequentially from either
of two Bernoulli processes.
The success probability of the first process is known, 
whereas that of the second is unknown. 
The second process is modeled as a Bayesian bandit whose state is the 
posterior distribution. 
They showed that the optimal policy is characterized by a break-even, critical
number,  which is a function of the number of remaining observations
and of the second process' state: one should continue sampling from
the second population as long as the current break-even value exceeds
the known success probability of the first process, and then
switch to the latter and keep sampling there until the $n$
observations are completed.
Such a break-even quantity is
the index of concern, although they did not use such a
term. Hence, their work introduced the  \emph{calibration
  approach} to bandit indexation, which has proven so fruitful in
later developments.

Jumping forward in time in this brief history of bandit indexation, 
\citet{whit88} significantly expanded the latter's scope beyond the
realm of classic bandits, by 
introducing an index for restless bandits ---
those that can change state while passive.
He did so by deploying a Lagrangian relaxation approach to the
intractable (cf.\ \citet{patsi99}) restless extension
of the classic multiarmed bandit problem under the average criterion. 
Whittle conjectured a form of asymptotic optimality of the resultant
index policy, which was established in \citet{webwei90, webwei91}
under certain conditions.

Yet, Whittle realized that
 existence of the index is not guaranteed for all restless
bandits: only for those that satisfy a so-called \emph{indexability} 
property. He stated in \citet{whit88}:
\begin{quote}
... one would very much like to have simple sufficient conditions for
indexability; at the moment, none are known.
\end{quote}

Such a state of affairs prompted the author to address that and other
issues on restless bandit indexation, as reported in the work surveyed
herein.
A cornerstone of our approach is the intuitive
concept introduced in \citet{nmmor06} of \emph{marginal productivity index} (MPI), which furnishes a
unifying framework for all the indices reviewed above as well as more
recent extensions.
Such a concept is grounded on insights drawn from the \emph{marginal
productivity theory} in economics developed at the end of the 19th
century by several researchers. See, e.g., the classic work by \citet{clark1899}.
The MPI of a project is a sound, intuitive priority index, as it measures the marginal 
productivity of work at each project state. 
As we will see, in the case of classic bandits the MPI reduces to an \emph{average
productivity} index. 
In a multi-project
setting, MPI policies  dynamically
assign higher priority 
to projects where work appears to be currently more productive. 
Besides being widely tractable and widely applicable, a growing body
of computational evidence
indicates that such index policies typically exhibit
a near-optimal performance and outperform conventional benchmark
policies derived from alternative approaches.

Several applications surveyed below are drawn from the domain of
optimal control of queueing systems. 
While the static optimization of such systems (cf.\ \citet{combeBoxma94} and
\citet{boxma95}) and some approaches to dynamic optimization have attracted
substantial research attention, emerging evidence suggests that,
within their scope,  
the MPI policies advocated in this paper can often yield
significant performance gains at a reduced
  computational expense.

The remainder of the paper is organized as follows.
Section \ref{s:rbi} reviews the key concepts and results of
the theory of restless bandit indexation, as well as of its
computational aspects as developed by the author extending Whittle's work. 
Section \ref{s:a} discusses applications to problems of admission control and
routing to parallel queues, scheduling a multiclass
make-to-order/make-to-stock queue, and scheduling a multiclass queue
with finite buffers.
Section \ref{s:rw} reports on more recent developments, involving theory, algorithms and
applications.
Finally, Section \ref{s:concl} concludes.

We remark that in the paper we use the terms ``bandit'' and
``project'' interchangeably.

\section{Restless bandit indexation: theory and computation}
\label{s:rbi}
We focus the following exposition on a discrete-time
single restless bandit model having a finite state space $N$, whose
one-period rewards and state-transition probabilities
under actions $a = 0$ (passive) and $a = 1$ (active)  at state $i$ are denoted 
by $R_i^a$ and $p_{ij}^a$, respectively. Rewards are discounted over
time with factor $0 < \beta < 1$.
The project is operated under a policy $\pi$, drawn from the class 
$\Pi$ of history-dependent randomized policies $\Pi$.
We denote by $X(t)$ and $a(t)$ the project
state and action processes,
respectively.

\subsection{Indexability and the MPI}
\label{s:iwimpi}
We evaluate a policy $\pi$ by means of two measures. The first is
the \emph{reward measure}
\[
f_{i_0}^\pi \triangleq \Ex_{i_0}^\pi\left[\sum_{t=0}^\infty R_{X(t)}^{a(t)} \beta^t\right],
\]
giving the expected total discounted value of rewards earned
over an infinite horizon starting at $i_0$. 
The second measure concerns the associated resource expenditure.
Thus, if $Q_j^a$ units of work are expended by taking action $a$ in
state $j$,  we use the \emph{work measure}
\[
g_{i_0}^\pi \triangleq \Ex_{i_0}^\pi\left[\sum_{t=0}^\infty Q_{X(t)}^{a(t)} \beta^t\right],
\]
giving the corresponding expected total discounted amount of work
expended.
We assume that such work-expenditure parameters satisfy
$Q_j^1 \geq Q_j^0 \geq 0$.

Note that such a setting allows the possibility that the two actions are 
identical in its resource consumption and dynamics at some states. We will find it convenient to identify the
states, if any,  where such is the case, 
\begin{equation}
\label{eq:uncst}
N^{\{0\}} \triangleq \big\{i \in N\colon Q_i^0 = Q_i^1 
\text{ and } p_{ij}^0 = p_{ij}^1, j \in N\big\},
\end{equation}
and call them \emph{uncontrollable}, while terming \emph{controllable}
the remaining states 
$N^{\{0, 1\}} \triangleq N \setminus N^{\{0\}}$.
We denote by $n \geq 1$ the number of controllable states, and 
adopt the convention that the passive action is taken at
uncontrollable states, which is reflected in the notation.

We will further refer to corresponding measures 
$f^\pi$ and $g^\pi$ obtained by drawing at random the initial state 
 according to an arbitrary
 positive probability mass function $p_{i} > 0$ for 
$i \in N$, i.e., $f^\pi \triangleq \sum_{i \in N} p_i f_i^\pi$ and 
$g^\pi \triangleq \sum_{i \in N} p_i g_i^\pi$.

Suppose that work is to be paid for at wage rate $\nu$, and
consider the  \emph{$\nu$-wage problem}
\begin{equation}
\label{eq:nuwp}
\max_{\pi \in \Pi} f^\pi - \nu g^\pi,
\end{equation}
which is to find an admissible project-operating policy that maximizes
the value of rewards earned minus labor costs incurred.
We will use (\ref{eq:nuwp}) as a 
 \emph{calibrating problem}, aimed
at measuring the marginal value of work at each project state.

Since (\ref{eq:nuwp}) is a finite-state and -action discounted MDP, 
standard results (cf.\ \citet{put94}) ensure existence of an optimal 
policy that is: (i) stationary deterministic; and (ii) independent of
the initial state. 
It is convenient to represent each such a policy by its \emph{active set} 
$S \subseteq N^{\{0, 1\}}$, which is the subset of states where it
prescribes to engage the project.
We will thus refer to the \emph{$S$-active policy} and write, e.g., 
$f^S$ and $g^S$. We can thus reduce (\ref{eq:nuwp}) to the 
\emph{combinatorial optimization problem} of finding an optimal active
set in the family of all subsets of $N^{\{0, 1\}}$, denoted by $2^{N^{\{0, 1\}}}$:
\begin{equation}
\label{eq:conuwp}
\max_{S \in 2^{N^{\{0, 1\}}}} f^S - \nu g^S.
\end{equation}

For every wage value $\nu$, the optimal policies are characterized by
the unique solution to the \emph{Bellman equations}
\begin{equation}
\label{eq:be}
\vartheta_i^*(\nu) = \max_{a \in \{0, 1\}} R_i^a - Q_i^a \nu + \beta \sum_{j
  \in N} p_{ij}^a \vartheta_j^*(\nu), \quad i \in N,
\end{equation}
where $\vartheta_i^*(\nu)$ denotes the optimal value of (\ref{eq:nuwp})
starting at $i$.
Hence, there exists a \emph{minimal optimal active set} $S^*(\nu)
\subseteq N^{\{0, 1\}}$ for (\ref{eq:nuwp}), which is characterized in
terms of (\ref{eq:be}) by 
\[
S^*(\nu) \triangleq \Big\{i \in N^{\{0, 1\}}\colon 
R_i^1 - Q_i^1 \nu + \beta \sum_{j
  \in N} p_{ij}^1 \vartheta_j^*(\nu) > R_i^0 - Q_i^0 \nu + \beta \sum_{j
  \in N} p_{ij}^0 \vartheta_j^*(\nu)\Big\}.
\]

Now, it appears reasonable that, at least in some models, active sets 
$S^*(\nu)$ should expand monotonically from the empty set $\emptyset$
to the full controllable state space $N^{\{0, 1\}}$ as the wage $\nu$
is decreased from $+\infty$ to $-\infty$. If such is the case, to each
controllable state $i$ will be attached 
a critical wage value $\nu_i^*$ below which $i$ enters 
$S^*(\nu)$.

\begin{definition}[Indexability; MPI]
\label{def:indmpi}
We say that the project is \emph{indexable} if there exists an 
\emph{index} $\nu_i^* \in \mathbb{R}$ for $i \in N^{\{0, 1\}}$ such that
\[
S^*(\nu) = \big\{i \in N^{\{0, 1\}}\colon \nu_i^* > \nu\big\}, \quad 
 \nu \in \mathbb{R}.
\]
In such a case we say that $\nu_i^*$ is the bandit's project.
\end{definition}

The concept of indexability was introduced by \citet{whit88}
 in the case $Q_i^a \equiv a$  under
 the long-run average criterion, in a formulation given in
 terms of optimal passive sets.
He also showed that, natural as it may seem, such a property should
not be taken for granted, as there are nonindexable projects.
The extension to the discounted criterion was carried out in
 \citet{nmaap01}.
The more general present setting was introduced in
 \citet{nmmp02}, where the index $\nu_i^*$ was first shown to measure
 the marginal value, or productivity, of work at each state, which,
 along with further extensions and results, 
 prompted our proposing the term MPI in 
\citet{nmmor06}.

\subsection{An achievable work-reward  view of indexability: geometric
 and economic insights}
\label{s:awrv}
\citet{nmmp02,nmmor06} introduces an achievable
work-reward approach to indexability that offers both geometric and
economic insights, having deep connections with multi-objective
optimization (cf.\ \citet{hl01}).

Consider the \emph{achievable work-reward performance region}
\[
\mathbb{H} \triangleq \big\{(g^\pi, f^\pi)\colon \pi \in \Pi\big\},
\]
which is the region spanned in the plane by work-reward performance
points under all admissible policies. 
Note that such a region is a convex polygon, given by the convex
hull 
of the finite collection of performance points $(g^S, f^S)$
achieved by stationary deterministic policies, as represented by their
active sets $S$:
\[
\mathbb{H} = \conv\Big(\big\{(g^S, f^S)\colon S \in 2^{N^{\{0, 1\}}}\big\}\Big).
\]

Of particular interest for our purposes is the \emph{upper boundary}
of such a region, 
\[
\bar{\partial} \mathbb{H} \triangleq \big\{(g, f) \in \mathbb{H}\colon
f^\pi \leq f \text{ for every } \pi \in \Pi\big\},
\]
as the indexability property concerns the latter's structure.
Thus, the project is indexable iff $\bar{\partial} \mathbb{H}$ is
characterized by a \emph{nested
active-set family}
\[
\mathcal{F}_0 \triangleq \big\{S_0, S_1, \ldots, S_n\big\},
\]
where $S_0 \triangleq \emptyset$, $S_n \triangleq N^{\{0, 1\}}$ and 
$S_k \triangleq \{i_1, \ldots, i_k\}$ for $1 \leq k \leq n$ satisfy 
\begin{equation}
\label{eq:gsklt}
g^{S_0} < g^{S_1} < \cdots < g^{S_n},
\end{equation}
and $i_1, \ldots, i_n$
is an ordering  of the project's $n$ controllable
states.

In such a case  the 
MPI has the evaluation
\begin{equation}
\label{eq:mpic}
\nu_{i_k}^* = \frac{f^{S_k} - f^{S_{k-1}}}{g^{S_k} - g^{S_{k-1}}},
\quad 1 \leq k \leq n.
\end{equation}
The latter representation shows that the MPI measures the  reward
vs.\ work trade-off rates or slopes in the upper boundary $\bar{\partial}
\mathbb{H}$, characterizing
 $\nu_{i_k}^*$ as the marginal value or
productivity of
work in the project at state $i_k$, which motivates our
using the term MPI.

We must emphasize that, 
while the measures $f^\pi$ and $g^\pi$ depend on the chosen
initial-state probabilities $p_i > 0$, the rates defining the MPI in 
(\ref{eq:mpic}) remain invariant under changes in the latter.
Further, whereas the present geometric approach requires in general
that such probabilities be positive, it can be shown that the MPI can
also be expressed as
\[
\nu_{i_k}^* = \frac{f_{i_k}^{S_k} - f_{i_k}^{S_{k-1}}}{g_{i_k}^{S_k} - g_{i_k}^{S_{k-1}}},
\quad 1 \leq k \leq n.
\]

The latter representation sheds light on the relation between the
index for restless and classic (nonrestless) bandits. Thus, for a
classic bandit with zero passive rewards we have
$f_{i_k}^{S_{k-1}} = g_{i_k}^{S_{k-1}} = 0$, and hence
\[
\nu_{i_k}^* = \frac{f_{i_k}^{S_k}}{g_{i_k}^{S_k}},
\quad 1 \leq k \leq n.
\]
Therefore, the MPI for a classic bandit reduces to an \emph{average 
productivity index}.

The above interpretation furnishes an intuitive economic justification
for use of MPI policies in a multi-project scenario. 
Thus, such policies seek to dynamically allocate work to those projects
that can make better use of it, using the MPI as a proxy --- as
it ignores interactions --- marginal productivity measure. 
Such a viewpoint draws and builds on the \emph{marginal
productivity theory} in economics, and its extensions which apply it
to optimal resource allocation. See,
e.g., \citet{clark1899}, \citet{koopm57} and \citet{kantor59}.

Notice further that the optimal value function $\vartheta_i^*(\nu)$ in (\ref{eq:be})
of an indexable project as above is given by
\begin{equation}
\label{eq:ovfrep}
\vartheta_i^*(\nu) = \max_{S \in \mathcal{F}_0} f_i^S - \nu g_i^S.
\end{equation}

\begin{figure}[htb]
\centering
\begin{psfrags}%
\psfragscanon%
%
\psfrag{s03}[l][l]{\color[rgb]{0,0,0}\setlength{\tabcolsep}{0pt}\begin{tabular}{l}$\emptyset$\end{tabular}}%
\psfrag{s04}[l][l]{\color[rgb]{0,0,0}\setlength{\tabcolsep}{0pt}\begin{tabular}{l}$\{1\}$\end{tabular}}%
\psfrag{s05}[l][l]{\color[rgb]{0,0,0}\setlength{\tabcolsep}{0pt}\begin{tabular}{l}$\{2\}$\end{tabular}}%
\psfrag{s06}[l][l]{\color[rgb]{0,0,0}\setlength{\tabcolsep}{0pt}\begin{tabular}{l}$\{3\}$\end{tabular}}%
\psfrag{s07}[l][l]{\color[rgb]{0,0,0}\setlength{\tabcolsep}{0pt}\begin{tabular}{l}$\{1, 2\}$\end{tabular}}%
\psfrag{s08}[l][l]{\color[rgb]{0,0,0}\setlength{\tabcolsep}{0pt}\begin{tabular}{l}$\{1, 3\}$\end{tabular}}%
\psfrag{s09}[l][l]{\color[rgb]{0,0,0}\setlength{\tabcolsep}{0pt}\begin{tabular}{l}$\{2, 3\}$\end{tabular}}%
\psfrag{s10}[l][l]{\color[rgb]{0,0,0}\setlength{\tabcolsep}{0pt}\begin{tabular}{l}$\{1, 2, 3\}$\end{tabular}}%
\psfrag{s11}[t][t]{\color[rgb]{0,0,0}\setlength{\tabcolsep}{0pt}\begin{tabular}{c}$g^{\pi}$\end{tabular}}%
\psfrag{s12}[b][b]{\color[rgb]{0,0,0}\setlength{\tabcolsep}{0pt}\begin{tabular}{c}$f^{\pi}$\end{tabular}}%
%
\psfrag{x01}[t][t]{0}%
\psfrag{x02}[t][t]{0.1}%
\psfrag{x03}[t][t]{0.2}%
\psfrag{x04}[t][t]{0.3}%
\psfrag{x05}[t][t]{0.4}%
\psfrag{x06}[t][t]{0.5}%
\psfrag{x07}[t][t]{0.6}%
\psfrag{x08}[t][t]{0.7}%
\psfrag{x09}[t][t]{0.8}%
\psfrag{x10}[t][t]{0.9}%
\psfrag{x11}[t][t]{1}%
%
\psfrag{v01}[r][r]{0}%
\psfrag{v02}[r][r]{0.1}%
\psfrag{v03}[r][r]{0.2}%
\psfrag{v04}[r][r]{0.3}%
\psfrag{v05}[r][r]{0.4}%
\psfrag{v06}[r][r]{0.5}%
\psfrag{v07}[r][r]{0.6}%
\psfrag{v08}[r][r]{0.7}%
\psfrag{v09}[r][r]{0.8}%
\psfrag{v10}[r][r]{0.9}%
\psfrag{v11}[r][r]{1}%
%
\includegraphics[height=2.2in,width=6.1in,keepaspectratio]{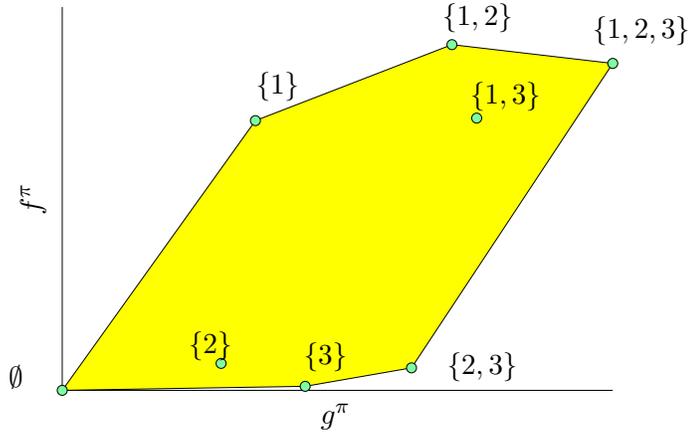}%
\end{psfrags}%
\caption{Achievable work-reward performance region of an indexable project.}
\label{fig:indpcl0}
\end{figure}

Two examples will help illustrate these ideas. 
Consider first the project with state space
$N = \{1, 2, 3\}$ and one-period work consumptions $Q_i^a = a$, 
discount factor 
$\beta = 0.9$, one-period active reward vector and one-period transition probabilities
\[
\vec{R}^1 =
\begin{bmatrix}
0.9016 \\ 0.10949 \\ 0.01055
\end{bmatrix}, 
\vec{P}^1 =
\begin{bmatrix}
0.2841 & 0.4827 & 0.2332 \\
0.5131 & 0.0212 & 0.4657 \\
0.4612 & 0.0081 & 0.5307
\end{bmatrix}, 
\vec{P}^0 =
\begin{bmatrix}
0.1810 & 0.4801 & 0.3389 \\
0.2676 & 0.2646 & 0.4678 \\
0.5304 & 0.2843 & 0.1853
\end{bmatrix},
\]
and one-period passive reward vector $\vec{R}^0 = \vec{0}$.
Figure \ref{fig:indpcl0} displays the achievable work-reward
performance region $\mathbb{H}$ for such an instance, where
points $(g^S, f^S)$ are labeled by their
  active sets $S$, and 
the initial-state distribution is uniform over $N$, i.e., $p_i =
1/3$ for $i \in N$.
The plot shows that this is an indexable project, relative to the nested
active-set family $\mathcal{F}_0 = \big\{\emptyset, \{1\},
\{1, 2\}, \{1, 2, 3\}\big\}$, which determines the region's upper
boundary $\bar{\partial}
\mathbb{H}$.
The MPI values of states $1$, $2$ and $3$ are given by
the successive trade-off vs.\ work rates or slopes in such an upper boundary:
\[
\nu_1^* = \frac{f^{\{1\}}-f^\emptyset}{g^{\{1\}}-g^\emptyset} > 
\nu_2^* = \frac{f^{\{1, 2\}}-f^{\{1\}}}{g^{\{1, 2\}}-g^{\{1\}}} > 
\nu_3^* = \frac{f^{\{1, 2, 3\}}-f^{\{1, 2\}}}{g^{\{1, 2, 3\}}-g^{\{1, 2\}}}.
\]

\begin{figure}[htb]
\centering
\begin{psfrags}%
\psfragscanon%
\psfrag{s03}[l][l]{\color[rgb]{0,0,0}\setlength{\tabcolsep}{0pt}\begin{tabular}{l}$\emptyset$\end{tabular}}%
\psfrag{s04}[l][l]{\color[rgb]{0,0,0}\setlength{\tabcolsep}{0pt}\begin{tabular}{l}$\{1\}$\end{tabular}}%
\psfrag{s05}[l][l]{\color[rgb]{0,0,0}\setlength{\tabcolsep}{0pt}\begin{tabular}{l}$\{2\}$\end{tabular}}%
\psfrag{s06}[l][l]{\color[rgb]{0,0,0}\setlength{\tabcolsep}{0pt}\begin{tabular}{l}$\{3\}$\end{tabular}}%
\psfrag{s07}[c][l]{\color[rgb]{0,0,0}\setlength{\tabcolsep}{0pt}\begin{tabular}{l}$\{1, 2\}$\end{tabular}}%
\psfrag{s08}[l][l]{\color[rgb]{0,0,0}\setlength{\tabcolsep}{0pt}\begin{tabular}{l}$\{1, 3\}$\end{tabular}}%
\psfrag{s09}[l][l]{\color[rgb]{0,0,0}\setlength{\tabcolsep}{0pt}\begin{tabular}{l}$\{2, 3\}$\end{tabular}}%
\psfrag{s10}[l][l]{\color[rgb]{0,0,0}\setlength{\tabcolsep}{0pt}\begin{tabular}{l}$\{1, 2, 3\}$\end{tabular}}%
\psfrag{s11}[t][t]{\color[rgb]{0,0,0}\setlength{\tabcolsep}{0pt}\begin{tabular}{c}$g^{\pi}$ \end{tabular}}%
\psfrag{s12}[b][b]{\color[rgb]{0,0,0}\setlength{\tabcolsep}{0pt}\begin{tabular}{c}$f^{\pi}$\end{tabular}}%
%
\psfrag{x01}[t][t]{0}%
\psfrag{x02}[t][t]{0.1}%
\psfrag{x03}[t][t]{0.2}%
\psfrag{x04}[t][t]{0.3}%
\psfrag{x05}[t][t]{0.4}%
\psfrag{x06}[t][t]{0.5}%
\psfrag{x07}[t][t]{0.6}%
\psfrag{x08}[t][t]{0.7}%
\psfrag{x09}[t][t]{0.8}%
\psfrag{x10}[t][t]{0.9}%
\psfrag{x11}[t][t]{1}%
%
\psfrag{v01}[r][r]{0}%
\psfrag{v02}[r][r]{0.1}%
\psfrag{v03}[r][r]{0.2}%
\psfrag{v04}[r][r]{0.3}%
\psfrag{v05}[r][r]{0.4}%
\psfrag{v06}[r][r]{0.5}%
\psfrag{v07}[r][r]{0.6}%
\psfrag{v08}[r][r]{0.7}%
\psfrag{v09}[r][r]{0.8}%
\psfrag{v10}[r][r]{0.9}%
\psfrag{v11}[r][r]{1}%
\includegraphics[height=2.2in,width=3in,keepaspectratio]{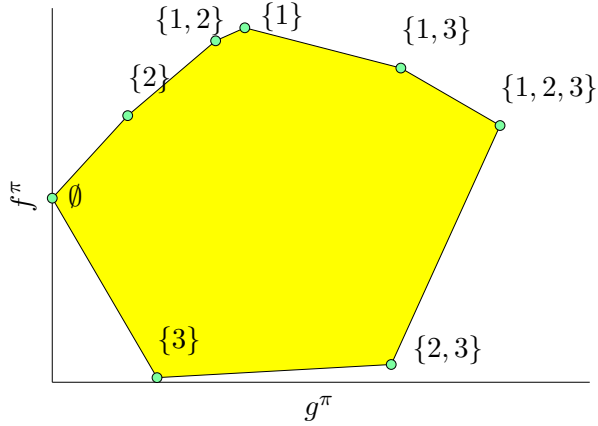}%
\end{psfrags}%
\caption{Achievable work-reward performance region of a nonindexable project.}
\label{fig:nonind1}
\end{figure}

Consider now the project instance having the same states,
initial-state distribution and discount factor, transition probabilities
\[
\vec{P}^1 =
\begin{bmatrix}
0.7796 & 0.0903 & 0.1301 \\
0.1903 & 0.1863 & 0.6234 \\
0.2901 & 0.3901 & 0.3198
\end{bmatrix}, \,
\vec{P}^0 =
\begin{bmatrix}
0.1902 & 0.4156 & 0.3942 \\
0.5676 & 0.4191 & 0.0133 \\
0.0191 & 0.1097 & 0.8712
\end{bmatrix},
\]
and reward vectors
\[
\vec{R}^1 = 
\begin{bmatrix} 
0.9631 \\ 0.7963 \\ 0.1057
\end{bmatrix}, \,
\vec{R}^0 = 
\begin{bmatrix} 
0.458 \\ 0.5308 \\ 0.6873
\end{bmatrix}.
\]
Figure \ref{fig:nonind1} displays the achievable work-reward
performance region for this instance. 
The plot reveals that this project is  nonindexable, since there is no
nested active-set family that determines the region's upper boundary.

\subsection{PCL-indexability conditions and adaptive-greedy index algorithm}
\label{s:sicaga}
While testing for indexability of a given restless bandit instance
is a conceptually simple task, as it can be solved, e.g.,  by 
visual inspection of 
plots such as those in Figures \ref{fig:indpcl0} and
\ref{fig:nonind1}, researchers will more often be interested 
in establishing \emph{analytically} that a particular  model 
arising in some
application is indexable under a suitable parameter range.
The latter task is, in contrast, generally far from trivial.
It would thus be useful to have tractable
 sufficient conditions for indexability that are widely applicable.
\citet{nmaap01,nmmp02,nmmor06} introduces, develops and deploys
the first such conditions, along with
a corresponding index algorithm, which we review next.

For a restless bandit as above, given an action $a \in \{0, 1\}$ and 
an active set $S \subseteq N^{\{0, 1\}}$, denote by $\langle a, S\rangle$ the
policy that takes action $a$ in the initial period and adopts the 
$S$-active policy thereafter. 
In addition to the work and reward measures discussed before, let us
now define
the \emph{marginal work measure}
\[
w_i^S \triangleq g_i^{\langle 1, S\rangle} - g_i^{\langle 0, S\rangle}
= Q_i^1 - Q_i^0 + \beta \sum_{j \in N} \big(p_{ij}^1 - p_{ij}^0\big) 
 g_j^S,
\]
and the \emph{marginal reward measure}
\[
r_i^S \triangleq f_i^{\langle 1, S\rangle} - f_i^{\langle 0, S\rangle}
= R_i^1 - R_i^0 + \beta \sum_{j \in N} \big(p_{ij}^1 - p_{ij}^0\big) 
 f_j^S.
\]
Notice that $w_i^S$ (resp.\ $r_i^S$) measures 
the marginal increment in work expended (resp.\ in value of rewards
earned)
that results from 
working instead of resting in the initial period starting at $i$, 
provided that the $S$-active policy is adopted afterwards.

Further, if $w_i^S \neq 0$, define the \emph{marginal productivity
  measure}
\[
\nu_i^S \triangleq \frac{r_i^S}{w_i^S}.
\]

Recall now the characterization of indexability discussed 
in Sect.\ \ref{s:awrv}. Typically, i.e., unless the state space is
linearly ordered, we will not be able to identify a priori the nested
active-set family $\mathcal{F}_0$ determining the achievable 
work-reward performance region's upper boundary.
Yet, often we can draw on intuition or experimentation to 
\emph{guess} the structure of optimal policies for the particular
model at hand, in the form of an active-set family $\mathcal{F} 
\subseteq N^{\{0, 1\}}$ that
\emph{contains} $\mathcal{F}_0$, i.e., $\mathcal{F} \supseteq
\mathcal{F}_0$,
 for a range of model parameters.
In fact, $\mathcal{F}$ will often be much larger than $\mathcal{F}_0$.
In the terminology of combinatorial optimization,
 $(N^{\{0, 1\}}, \mathcal{F})$ is a \emph{set system} on \emph{ground set}
$N^{\{0, 1\}}$ having $\mathcal{F}$ as its family of \emph{feasible
 sets}.

Algorithmic considerations, namely the requirement that we can build
 our way up from the empty set 
towards a given set $S \in \mathcal{F}$ through successive 
single-state augmentations, as well as the symmetric requirement for 
reaching $S$ through successive single-state removals from 
$N^{\{0, 1\}}$,  lead us to impose some natural conditions
 on such a set system, for which we need the following concepts.
For $S \in \mathcal{F}$, define the \emph{inner boundary of $S$
 relative to $\mathcal{F}$} by 
\[
\partial^{\textup{in}}_{\mathcal{F}} S \triangleq 
\big\{i \in S\colon S \setminus \{i\} \in \mathcal{F}\big\}.
\]
Define further the \emph{outer boundary of $S$
 relative to $\mathcal{F}$} by 
\[
\partial^{\textup{out}}_{\mathcal{F}} S \triangleq 
\big\{i \in N^{\{0, 1\}} \setminus S\colon S \cup \{i\} \in \mathcal{F}\big\}.
\]

\begin{assumption}
\label{ass:nfreq}
Set system $(N^{\{0, 1\}}, \mathcal{F})$ satisfies the following
conditions: 
\begin{itemize}
\item[(i)] $\emptyset, N^{\{0, 1\}} \in \mathcal{F}$;
\item[(ii)] for $\emptyset \neq S \in \mathcal{F}$,
  $\partial^{\textup{out}}_{\mathcal{F}} S \neq \emptyset$;
\item[(iii)] for $N^{\{0, 1\}} \neq S \in \mathcal{F}$,
  $\partial^{\textup{in}}_{\mathcal{F}} S \neq \emptyset$.
\end{itemize}
\end{assumption}

Consider now the \emph{adaptive-greedy algorithm}
$\mathrm{AG}_{\mathcal{F}}$ shown in Table \ref{fig:agaf}.
In essence, referring to the geometric viewpoint in Sect.\
\ref{s:awrv},
this algorithm seeks to traverse the upper boundary of the 
achievable work-reward performance region, building up the successive active sets $S_k$ forming the nested family 
$\mathcal{F}_0$ that determines such a boundary. To do so it restricts
attention to active sets drawn from the given family $\mathcal{F}$. 
Notice further that the algorithm aims to traverse such a boundary
from left to right, so that the successive index values or slopes in 
such a frontier are computed in nonincreasing order, i.e., it is a 
\emph{top-down} index algorithm.
Also, the algorithm is only well defined when the computed 
marginal productivity rates have nonzero denominators.
The output consists of an ordered string $i_1, \ldots, i_n$
of the $n$ controllable states in 
$N^{\{0, 1\}}$, along with corresponding
index values $\nu^*_{i_k}$. 

\begin{table}[htb]
\caption{Adaptive-greedy index algorithm 
$\mathrm{AG}_{\mathcal{F}}$.}
\begin{center}
\fbox{%
\begin{minipage}{1.8in}
\textbf{ALGORITHM} $\mathrm{AG}_{\mathcal{F}}$: \\
\textbf{Output:}
$\{i_k, \nu^*_{i_k}\}_{k=1}^{n}$
\begin{tabbing}
$S_0 := \emptyset$ \\
\textbf{for} \= $k := 1$ \textbf{to}  $n$ \textbf{do} \\
\> \textbf{pick}  
 $i_k \in \argmax
      \big\{\nu^{S_{k-1}}_{i}\colon
                 i \in \partial^{\textup{out}}_{\mathcal{F}} S_{k-1}\big\}$ \\
 \>  $\nu^*_{i_k} := 
 \nu^{S_{k-1}}_{i_k}$;  \, $S_{k} := S_{k-1} \cup \{i_k\}$ \\
\textbf{end} \{ for \}
\end{tabbing}
\end{minipage}}
\end{center}
\label{fig:agaf}
\end{table}

We next use such an algorithm to define a certain class of restless
bandits. Note that the acronym ``PCL'' refers to the 
\emph{partial conservation laws} introduced in \citet{nmaap01}.
\begin{definition}[PCL$(\mathcal{F})$-indexability]
\label{def:pclib}
We say that a bandit is \emph{PCL$(\mathcal{F})$-indexable} if it
satisfies the following conditions:
\begin{itemize}
\item[(i)] Positive marginal work: $w_i^S > 0$ for $i \in N^{\{0,
    1\}}, S \in \mathcal{F}$; 
\item[(ii)] Monotone nonincreasing index computation: the index values produced by 
algorithm $\mathrm{AG}_{\mathcal{F}}$
  satisfy
\[
\nu^*_{i_1} \geq \nu^*_{i_2} \geq \cdots \geq \nu^*_{i_n}.
\]
\end{itemize}
\end{definition}

Note that part (i) of Definition \ref{def:pclib} ensures, along with 
Assumption \ref{ass:nfreq}, that the algorithm is well defined.
The interest of the class of PCL$(\mathcal{F})$-indexable bandits 
is based on the following result, proven in
\citet[Cor.\ 2]{nmaap01},
\citet[Th.\ 6.3]{nmmp02} and \citet[Th.\ 4.1]{nmmor06} in increasingly general settings.

\begin{theorem}
\label{the:pclii}
A PCL$(\mathcal{F})$-indexable bandit is indexable and 
algorithm $\mathrm{AG}_{\mathcal{F}}$ gives its MPI.
\end{theorem}

From the point of view of the practical application of Theorem
\ref{the:pclii} to a particular
model, one would first set out to establish analytically satisfaction 
of Definition \ref{def:pclib}(i). 
Note that it might well happen that $w_i^S < 0$ for some active sets
not in $\mathcal{F}$. Yet, as stated, it suffices to prove positivity of
marginal work measures for active sets drawn from family
$\mathcal{F}$.
Then, one would check for satisfaction of Definition
\ref{def:pclib}(ii).
This can be done either computationally, simply by running the
algorithm and testing whether the index is computed in nonincreasing
order, or analytically. We refer the reader, e.g., 
 to \citet{nmmp02,nmmor06,nmqs06} for examples of detailed analyses of 
specific models.

We next comment on an approach we have found useful to establish analytically
condition (ii) in Definition \ref{def:pclib}, once part (i) has been
proven. 
Thus, suppose we want to show that the $k$th and $(k+1)$th computed index values satisfy 
$\nu_{i_k}^* \geq \nu_{i_{k+1}}^*$ for any $1
  \leq k < n$.
Now, letting $S_{k-1}$ and $S_k =
  S_{k-1} \cup \{i_k\}$ be as in Table \ref{fig:agaf},
we can use \citet[Prop.\ 6.4(c)]{nmmp02} to write
\[
r_{i}^{S_k} - r_{i}^{S_{k-1}} =
\frac{r_{i_k}^{S_k}}{w_{i_k}^{S_k}} \big(w_{i}^{S_k} -
w_{i}^{S_{k-1}}\big), \quad i \in N^{\{0, 1\}},
\]
which is immediately reformulated using 
$\nu_{i_k}^* = \nu_{i_k}^{S_{k-1}}$ as
\[
\nu_{i}^{S_k} = \nu_{i_k}^* - 
\frac{w_{i}^{S_{k-1}}}{w_{i}^{S_{k}}} \big(\nu_{i_k}^* -
\nu_{i}^{S_{k-1}}\big), \quad i \in N^{\{0, 1\}}.
\]
From the latter identity along with condition (i)
we obtain
\begin{equation}
\label{eq:equivineq}
\nu_{i_k}^* \geq \nu_{i_{k+1}}^*
\Longleftrightarrow 
\nu_{i_k}^* \geq \nu_{i_{k+1}}^{S_{k-1}},
\end{equation}
which shows that it suffices to prove that $\nu_{i_k}^* \geq 
\nu_{i_{k+1}}^{S_{k-1}}$.

To illustrate, consider the case of a classic bandit --- where 
$N = N^{\{0, 1\}}$. Taking $\mathcal{F} \triangleq 2^N$ it is easily
shown (cf.\ \citet[Cor.\ 4]{nmaap01}) that Definition
\ref{def:pclib}(i) holds. To prove condition (ii) note that, by
construction, 
\[
\nu_{i_k}^* = \max \big\{\nu_i^{S_{k-1}}\colon i \in N \setminus
S_{k-1}\big\} \geq \nu_{i_{k+1}}^{S_{k-1}},
\]
and hence, by (\ref{eq:equivineq}), we obtain $\nu_{i_k}^* \geq 
\nu_{i_{k+1}}^*$.
Therefore, classic bandits are PCL$(2^N)$-indexable, or, in the
terminology used in \citet{nmaap01}, \emph{GCL-indexable}, as in such
a case they satisfy the \emph{generalized conservation laws} (GCL) in 
\citet{beni}. It is further shown in \citet[Cor.\ 5]{nmaap01} that
any restless bandit is GCL-indexable for small enough values of the discount factor.

The reader may wonder about the intuitive interpretation of
Definition \ref{def:pclib}(i), besides that suggested by definition of
the $w_i^S$'s.
Further insights into such an issue are given in
 \citet[Prop.\ 6.2]{nmmp02}, where it is shown
that the condition can be reformulated in terms of work
measures as follows. For $S \in \mathcal{F}$, 
\begin{equation}
\label{eq:equivrefwisp}
\begin{split}
g_i^S & < g_i^{S \cup \{i\}}, \quad i \in N^{\{0, 1\}} \setminus S \\
g_i^S & > g_i^{S \setminus \{i\}}, \quad i \in S.
\end{split}
\end{equation}
Note that (\ref{eq:equivrefwisp}) represents a
regularity property of 
work measures for active sets in $\mathcal{F}$
whereby, starting at a state $i$: augmenting an active set in 
$\mathcal{F}$ by adding
$i$ leads to an increase in work expended, whereas shrinking an active set
by removing $i$ leads to a decrease in work expended.

Another insightful representation of the index is given in \citet[Th.\
  6.4]{nmmp02} and \citet[Th.\ 4.2]{nmmor06}. Letting 
$\mathcal{F}_0 = \{S_0, S_1, \ldots, S_n\}$ be the active-set family 
produced by the algorithm, it holds that, for $1 \leq k \leq n$,
\begin{equation}
\label{eq:insr1}
\max_{j \in N^{\{0, 1\}} \setminus S_{k-1}}
\nu_{j}^{S_{k-1}} = 
\nu_{i_k}^{S_{k-1}}
= \nu_{i_k}^* = \nu_{i_k}^{S_{k}} = \min_{j \in S_k} \nu_{j}^{S_{k}}
\end{equation}
which is equivalent to the more intuitive reformulation
\begin{equation}
\label{eq:insr2}
\max_{j \in N^{\{0, 1\}} \setminus S_{k-1}}
\frac{f^{S_{k-1} \cup \{j\}} - f^{S_{k-1}}}{g^{S_{k-1} \cup \{j\}} -
  g^{S_{k-1}}} =
\frac{f^{S_{k}} - f^{S_{k-1}}}{g^{S_{k}} -
  g^{S_{k-1}}}  = \nu_{i_k}^* = 
\min_{j \in S_k} 
\frac{f^{S_k} - f^{S_k \setminus \{j\}}}{g^{S_k} - g^{S_k \setminus \{j\}}}.
\end{equation}
Such relations 
characterize the index as a \emph{locally optimal} marginal
productivity rate. Note that (\ref{eq:insr2}) has a clear geometric 
interpretation in the setting of the achievable work-reward region
approach in Sect.\ \ref{s:awrv}.

Further, we have found that, in some models, marginal work measures
satisfy the following monotonicity condition:
\begin{equation}
\label{eq:mc1}
\text{if } S, S' \in \mathcal{F} \text{ and } i \in S \subset S'
\text{ then } w_i^S \leq w_i^{S'}.
\end{equation}
Under (\ref{eq:mc1}), it is shown in \citet[Th.\ 4.7]{nmmp02} that
the index has the alternative representation
\begin{equation}
\label{eq:indr1}
\nu_i^* = \max_{S \in \mathcal{F}_0\colon i \in S} \nu_i^S = \max_{S
  \in \mathcal{F}_0\colon i \in S} \frac{f_i^S - f_i^{S \setminus
  \{i\}}}{g_i^S - g_i^{S \setminus \{i\}}},
\end{equation}
which is closely related to the Gittins index representation given in 
\citet{gi79} as an optimal average reward rate relative to stopping times. 
Thus, in the special case of a classic bandit with zero passive rewards, (\ref{eq:indr1})
reduces to 
\[
\nu_i^* = \max_{S
  \in \mathcal{F}_0\colon i \in S} \frac{f_i^S}{g_i^S},
\]
whereas the Gittins index representation referred to above appears,
when stationary deterministic stopping times are formulated in terms of their active sets, as
\[
\nu_i^* = \max_{S \in 2^N\colon i \in S} \frac{f_i^S}{g_i^S}.
\]

Such a result is extended in \citet[Th.\ 4.3]{nmmor06},
where it is shown that, if marginal work measures are
\emph{wedge-shaped}, meaning that they satisfy (\ref{eq:mc1}) and the
condition
\begin{equation}
\label{eq:mc2}
\text{if } S, S' \in \mathcal{F}, S \subset S' \text{ and }  i \not\in
S' 
\text{ then } w_i^S \geq w_i^{S'},
\end{equation}
then the index also has the representation
\begin{equation}
\label{eq:indr2}
\nu_i^* = \min_{S \in \mathcal{F}_0\colon i \not\in S} \nu_i^S = \min_{S
  \in \mathcal{F}_0\colon i \not\in S} \frac{f_i^{S \cup \{i\}} - f_i^{S}}{g_i^{S \cup \{i\}} - g_i^{S }}.
\end{equation}

It is worth outlining the evolution of ideas that led to our
introduction of PCLs and algorithm $\mathrm{AG}_{\mathcal{F}}$.
The early roots are in the pioneering work of \citet{kl}, who introduced 
an adaptive-greedy algorithm based on LP duality to compute the
average-optimal index policy for scheduling a multiclass $M/G/1$ queue with 
Bernoulli feedback.
\citet{vawabu} extended the scope of the 
adaptive-greedy algorithm to compute the Gittins index.
Then, drawing on earlier research on work conservation laws in multiclass queueing
systems by \citet{kleinr65}, \citet{coffmit80}, \citet{fedgro88},
\citet{shaya92} and \citet{tsouc91}, \citet{beni} developed a general
polyhedral framework 
for the study of indexation encompassing such previous work, based on the
unifying concept of GCL mentioned above.
In short, a generic scheduling problem obeying such laws is solved optimally
by an index policy, where the index is computed by an adaptive-greedy
algorithm. In essence, the latter correspond to the case of 
$\mathrm{AG}_{\mathcal{F}}$ where $\mathcal{F}$ consists of all
subsets of the ground set.
Yet, unlike classic bandits, restless bandits
do not necessarily satisfy GCLs.
This prompted the author to develop a
framework that applies to the
latter, based on a priori identification of a suitable active-set family 
$\mathcal{F}$. The resultant framework, termed 
PCL$(\mathcal{F})$-indexability after its grounding on 
satisfaction of partial conservation laws, was introduced
in \citet{nmaap01} along with algorithm $\mathrm{AG}_{\mathcal{F}}$
and Theorem \ref{the:pclii}.
To be precise, that paper gave the first version of algorithm
$\mathrm{AG}_{\mathcal{F}}$, yet in a different formulation than that
in Table \ref{fig:agaf}. The present, equivalent formulation was given
in \citet{nmmp02}, in an extended setting that grounds the approach on 
polyhedral LP methods. 

We must remark that, as shown in Table \ref{fig:agaf},
$\mathrm{AG}_{\mathcal{F}}$ is not really an algorithm, but an 
\emph{algorithmic scheme}, as no implementation details are given.
For a discussion of actual implementations see
Sect.\ \ref{s:mpicfic}, which reviews results from \citet{nmbfsind07}.
The latter paper reveals and exploits the deep connection between the
adaptive-greedy index algorithm and the classic parametric-objective 
simplex method of \citet{gassaaty55}.

\subsection{Extensions}
\label{s:se}

\subsubsection{Indexation in a general framework and 
  the law of diminishing marginal returns}
\label{s:gfirldmr}
The indexation theory for restless bandits 
is developed in a general framework in
terms of generic work and reward (or cost) measures in
\citet[Sect.\ 3--4]{nmmor06}, which highlights fundamental properties while hiding
ancillary model-specific details.
Such an approach reveals the fundamental connection between the
indexability property and the classic economic law of
\emph{diminishing marginal returns}, or \emph{diminishing marginal
  productivity}, whereby, as usage of a resource increases, its 
marginal productivity diminishes. 
Thus, indexable projects are precisely those that obey such a law as
it applies to the work expended on them, in such a way that there is 
a well defined marginal value of work at every state. See Th.\ 3.1 in
that paper.

\subsubsection{Semi-Markov bandits and bandits with a countable state space}
\label{s:smrb}
\citet{nmmor06} also shows how to adapt the theory to semi-Markov
restless bandits, which is relatively straightforward as the latter
are readily reformulated into a discrete-stage setting by standard methods. 
Further, that paper extends the indexation theory to bandits with a
countable state space,  yet restricting attention to the case where the
state space is linearly ordered.
The analysis of the general countable-state case raises considerable
technical difficulties, as, e.g., it is not clear how to ensure that
the state sequence produced by the adaptive-greedy algorithm  traverses the entire controllable state space.

\subsubsection{Average, discounted, bias and mixed criteria}
\label{s:aoc}
The theory of restless bandit indexation was introduced by 
\citet{whit88} under the 
average criterion, and then extended to the discounted criterion 
in \citet{nmaap01}.
Such criteria fit as special cases in the general indexation theory 
discussed in \citet[Sect.\ 3--4]{nmmor06}. 
We have found that the resultant flexibility furnishes an expanded scope that is relevant in applications. 

Consider, e.g., the classic problem of scheduling a multiclass $M/G/1$ 
to minimize average linear holding costs, which is well known to be
solved by the $c \mu$ index rule (cf.\ \citet{coxsmt}). 
The problem is immediately formulated as a restless bandit problem by viewing each class' queue as a restless
bandit project. 
Yet, as noticed in \citet[Chapt.\ 14.7]{whit96}, and also in
\citet{vewe} in a related  multiclass make-to-stock model,
such bandits are not indexable under the average criterion. As the
latter authors put it: 
\begin{quote}
In contrast, the backorder problem is not indexable. $\nu(x)$ does not exist (i.e., equals $-\infty$) for all $x$. The difficulty is
that $\nu$ is a Lagrange multiplier for the constraint on the time-average number of active arms. For the backorder problem,
any stable policy must serve a time-average of $\rho$ classes, so
relaxing this constraint does not change the optimal value, 
and the Lagrange multiplier does not exist. In fact, no scheduling problem with a fixed utilization will be indexable.
\end{quote}

We first showed in \citet{nmaller03}, in the setting of a multiclass
make-to-order/make-to-stock queue with convex stock and backlog
holding cost rates, that such a difficulty is
resolved by defining indexability relative to a \emph{mixed
  criteria} version of
$\nu$-wage problem (\ref{eq:nuwp}), where reward measure $f_i^\pi$ is
evaluated by the average criterion while work measure $g_i^\pi$ is
evaluated by the \emph{bias} criterion of \citet{blackwell62}.
For a bandit representing an  $M/G/1$ queue with utilization
factor $\rho < 1$ subject to service control, the latter is defined by 
\begin{equation}
\label{eq:biasc}
g_i^\pi \triangleq \Ex_i^\pi\left[\int_0^{\infty} \big\{a(t) -
  \rho\big\} \, dt\right],
\end{equation}
and thus gives the expected total cumulative excess work expended over
the nominal allocation.

These ideas are developed in \citet[Sect.\ 5.3]{nmmor06}.
Such an approach reveals a new ground of practical application of bias
optimality and mixed criteria in MDPs, which previously had 
been mostly considered topics of theoretical interest.
See \citet{lewput02}, and \citet{feinsh02}.

The scope of indexability is further extended to a pure
bias criterion in \citet[Sect.\ 2.3 and 4.2]{nmqs06}, motivated by
the analysis of a model involving the dynamic scheduling of a multiclass
queue with finite buffers, which is nonindexable under the average criteria.
Of course, the indices defined under the mixed average-bias and the pure
bias criteria are related to their discounted counterparts through a
vanishing-discount approach.

\section{Applications}
\label{s:a}
   \subsection{Control of admission and routing to parallel queues}
   \label{s:carpq}
In \citet{nmmp02}, the scope of indexation
is extended by allowing work-consumption rates to be
state-dependent. Besides introducing new LP-based polyhedral methods
for the analysis and computation of the resultant indices, the
approach is illustrated to develop new index policies for admission
control and routing to parallel queues.
The individual projects of concern are Markovian birth-death queues subject to
control of admission, which have
finite buffer space,
state-dependent arrival and service rates, and which incur nonlinear holding costs as well as
rejection costs. 
While such problems had been the subject of extensive research
typically based on DP analyses, 
as surveyed in \citet{stidh85}, the indexation approach yields new
insights, along with a fruitful connection with routing
problems. 

The single-queue admission control model addressed in
that paper is an extension of that in \citet{chenyao}. 
It incorporates arrival rates $\lambda_i$, service rates 
$\mu_i$, and holding cost rates $h_i$ that depend on the queue's
current state $i$, taken as the number of customers in system.
We show that, if the difference $\mu_i - \lambda_i$ is concave
nondecreasing and $h_i$ is convex nondecreasing in $i$, then the
model is indexable both under the discounted and the average
criterion. 
Yet, such an indexability result holds relative to an unconventional work
measure $g_i^\pi$ introduced in \citet{nmmp02}, which gives the expected total discounted
number or the expected long-run average rate per unit time, as the case may be, of arriving customers that are rejected due either to their finding
a full buffer upon arrival, or to their finding a \emph{closed entry
  gate} --- imagine that admission control is implemented through a
gatekeeper that opens or shuts an entry gate to the system. 
Notice that, when the buffer is full, opening or closing the entry gate,
i.e., taking the passive or the active action, have identical
consequences in terms of work expended and state dynamics. 
The corresponding state is hence uncontrollable (cf.\ (\ref{eq:uncst})).

Under the stated conditions, the model is shown to be PCL-indexable,
and hence 
indexable, with an index $\nu_i^*$ that is nondecreasing in $i$.
Such an index characterizes the optimal policies to the single-queue admission
control problem that incorporates a cost $\nu$ per customer rejected:
it is optimal to reject a customer who arrives in state $i < n$, where
$n$ is the buffer size,  iff
$\nu_i^* \geq \nu$. Notice that optimal policies are hence of
\emph{threshold} type: reject arrivals if the number of customers in system $i$ is
large enough. 
While the conventional approach to such problems focused on explicit 
determination of such thresholds, the indexation approach determines
them implicitly.

To compute the index, an efficient upwards recursion is given that generates
index values starting at $\nu_0^*$. This shows that the index does not
depend on the buffer size, and hence extends to a corresponding model with infinite
buffer space.

Note further that the index is a state-dependent measure of the value of
rejecting an arrival, which suggests a natural index policy in the
corresponding multi-project model, involving the 
admission control and routing to multiple queues in parallel.
Thus, consider now a model where customers arrive as a Poisson stream
with rate $\lambda$, and each rejected customer incurs a cost $\nu$.
The controller decides whether or not to reject each customer upon
arrival. If admitted, the customer is to be routed irrevocably to one
of $K$ queues in parallel, among those that are not full. 
Each such queue is modeled as in the single-queue model discussed
above, yet now assuming a constant arrival rate equal to $\lambda$.
Hence, if all the individual queues' buffers are full, an arriving customer is
necessarily rejected.
If we now denote by $\nu_k^*(i_k)$ the index for queue $k$ in state
$i_k$, the resultant admission control and routing policy is as
follows. 
Upon an arrival that finds each queue $k$ in state $i_k$:
reject the customer either if all buffers are full or if 
$\nu_k^*(i_k) \geq \nu$ for each nonfull queue $k$; otherwise, route the
customer to a nonfull queue of \emph{minimum index}
$\nu_k^*(i_k)$.

One may also use the index policy in the problem
version where only the
routing control capability is enabled, and  in a model where
queues have
infinite buffer space.
The proposed index policies reduce to the 
\emph{shortest (nonfull) queue routing} policy in the symmetric cases
in which this is optimal. 
See \citet{winston77}, \citet{johri89}, and \citet{hk90}.
The index policy also recovers the optimal policy in the
nonsymmetric routing model solved in \citet{dermanetal80}, where each queue
has a single buffer space.

In the simplest case of an $M/M/1$ queue with parameters $\lambda$ and
$\mu$  under
linear holding costs $h_i \triangleq c i$, the index is
shown in \citet[Sect.\ 7]{nmmp02} to have the evaluation
\begin{equation}
\label{eq:nukplhc}
\nu_j^* = 
\frac{c}{\mu} \sum_{i=1}^{j+1} \left(1 + \cdots +
  \rho^{i-1}\right) =
\begin{cases}
\displaystyle \frac{c}{\mu} \, \left[\frac{\rho^{j+2}-1}{(\rho-1)^2} -
 \frac{j+2}{\rho-1} 
 \right]
 & \text{if } \rho \neq 1 \\ \\
\displaystyle \frac{c}{\mu} \, \frac{(j+1) \, (j+2)}{2} & \text{if } \rho = 1,
\end{cases}
\end{equation}
where $\rho \triangleq \lambda/\mu$ is the utilization factor. Note
that the latter will typically be greater than one for each queue in a routing model
where there are multiple queues in parallel. A closed formula is also
given in that paper for the quadratic holding cost case $h_i
\triangleq c i^2$.
Further, in  \citet[Remark 7.3]{nmmp02} it is noted that, under the
average criterion, such a model is PCL-indexable provided only that
the holding cost rate is nondecreasing, i.e., the convexity assumption
can be then dropped.

In \citet{nmnetcoop07} we have investigated  the
applicability of such index policies to a model
involving the dynamic control of admission and routing to parallel
multi-server loss queues with reneging. 
The policies derived from 
restless bandit indexation are new. More importantly, they appear to be useful, as
the preliminary experiments reported in \citet{nmnetcoop07} reveal both a near-optimal
performance and substantial gains against conventional benchmark
policies
 on the instances investigated.

We remark that a similar model had been addressed in
\citet{kallmescass95} via a static optimization approach (cf.\
\citet{combeBoxma94} and \citet{combeBoxma94}), which is
hindered by the lack of closed formulae for the functions to be
optimized.
Other approach that has attracted substantial research attention is
based on improving upon the optimal static allocation by carrying out 
one step of the policy iteration algorithm. 
\citet{kr88, kr90} deploys such an approach to problems of optimal
routing to parallel queues, obtaining dynamic index policies. 

As part of a paper currently under preparation, we are conducting
extensive computational experiments comparing the MPI policies both
against optimal and alternative benchmark policies. 
We advance next an illustrative result. 
Consider the classic problem of routing a Poisson stream of customers
arriving with rate $\lambda$ to two $M/M/1$ queues in parallel, with
queue $k$'s server having service rate $\mu_k$ for $k = 1, 2$.
The aim is to find a routing policy that minimizes the  long-run
average customer
sojourn time.
To investigate computationally the performance of alternative policies
for such a system, we use a modified system where each queue has a finite
buffer capable of holding up to $30$ jobs, waiting or in service.
Also, we fix the arrival rate to the value $\lambda = 1$, and let
service rates $(\mu_1, \mu_2)$ vary over the range $1/2 < \rho < 1$,
where $\rho = \lambda/(\mu_1 + \mu_2)$. Note that $\rho \approx 1$
corresponds to the heavy-traffic regime.
The parameters $\mu_k$ are varied over a finite grid of width $1/100$.

We thus obtain that, over such a region, the relative suboptimality
gap of the MPI routing policy reaches a maximum value of about
$1.55\%$, being substantially smaller in most of the region.
Also, the MPI policy always outperforms the classic \emph{join the
  shortest queue} routing rule, achieving relative performance gains
of over $84\%$ when both service rates are very different in
magnitude.
As for the \emph{individually optimal routing} rule, we find that, 
though the MPI
policy can be worse than it over a small parameter region, it can lose
no more than $0.37\%$ in relative performance. Yet, over most of the
parameter region the MPI policy performs better, attaining relative
performance gains of up to over $15\%$.

\begin{figure}[htb]
\centering
\begin{psfrags}%
\psfragscanon%
\psfrag{s05}[t][t]{\color[rgb]{0,0,0}\setlength{\tabcolsep}{0pt}\begin{tabular}{c}$\mu_1$\end{tabular}}%
\psfrag{s06}[b][b]{\color[rgb]{0,0,0}\setlength{\tabcolsep}{0pt}\begin{tabular}{c}$\mu_2$\end{tabular}}%
\psfrag{s08}[b][b]{\color[rgb]{0,0,0}\setlength{\tabcolsep}{0pt}\begin{tabular}{c}$\Delta^{\textup{MPI}, \textup{OSI}}$\end{tabular}}%
%
\psfrag{x01}[t][t]{0}%
\psfrag{x02}[t][t]{0.1}%
\psfrag{x03}[t][t]{0.2}%
\psfrag{x04}[t][t]{0.3}%
\psfrag{x05}[t][t]{0.4}%
\psfrag{x06}[t][t]{0.5}%
\psfrag{x07}[t][t]{0.6}%
\psfrag{x08}[t][t]{0.7}%
\psfrag{x09}[t][t]{0.8}%
\psfrag{x10}[t][t]{0.9}%
\psfrag{x11}[t][t]{1}%
\psfrag{x12}[t][t]{0}%
\psfrag{x13}[t][t]{0.2}%
\psfrag{x14}[t][t]{0.4}%
\psfrag{x15}[t][t]{0.6}%
\psfrag{x16}[t][t]{0.8}%
\psfrag{x17}[t][t]{1}%
\psfrag{x18}[t][t]{0.01}%
\psfrag{x19}[t][t]{1}%
\psfrag{x20}[t][t]{2}%
%
\psfrag{v01}[r][r]{0}%
\psfrag{v02}[r][r]{0.2}%
\psfrag{v03}[r][r]{0.4}%
\psfrag{v04}[r][r]{0.6}%
\psfrag{v05}[r][r]{0.8}%
\psfrag{v06}[r][r]{1}%
\psfrag{v07}[l][l]{$-0.68\%$}%
\psfrag{v08}[l][l]{$-0.01\%$}%
\psfrag{v09}[l][l]{$0.01\%$}%
\psfrag{v10}[l][l]{$1\%$}%
\psfrag{v11}[l][l]{$3\%$}%
\psfrag{v12}[l][l]{$5\%$}%
\psfrag{v13}[l][l]{$7.64\%$}%
\psfrag{v14}[r][r]{0.01}%
\psfrag{v15}[r][r]{1}%
\psfrag{v16}[r][r]{2}%
%
\includegraphics[height=2.2in,width=6.1in,keepaspectratio]{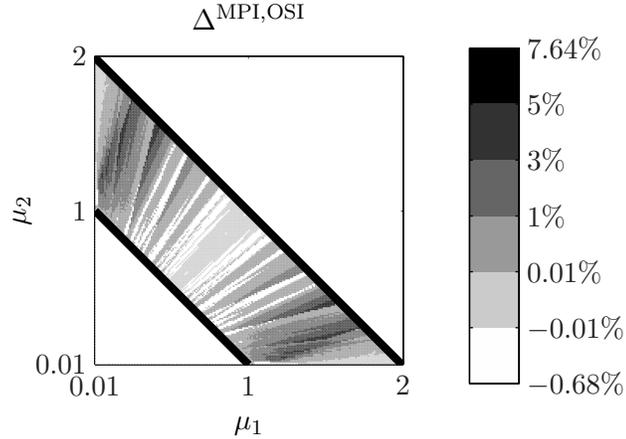}%
\end{psfrags}%
\caption{Relative performance gains of MPI over  OBS and OSI policies.}
\label{fig:picexp1_3}
\end{figure}

We have also compared the performance of the MPI routing policy
against the the policy obtained from the best static policy, i.e., 
\emph{optimal Bernoulli splitting} (OBS), 
 via \emph{one-step
policy improvement} (OSI), as proposed in \cite{kr90}.
The results are displayed in Figure \ref{fig:picexp1_3}.
We emphasize that in such results the OBS policy's objective was
exactly computed for an infinite-buffer system using closed formulae
as in \cite{kr90},
while the OSI and the MPI policies'
objectives were numerically computed for the finite-buffer approximation
mentioned above.
The figure shows a density plot of $\Delta^{\textup{MPI},
  \textup{OSIS}}$, the 
relative performance gain 
of the MPI over the OSI policy, within the parameter region of concern
delimited by the displayed thick lines corresponding to $\rho = 1/2$ and
$\rho = 1$. 
We find that in some instances
 the MPI policy's performance can be worse than the OSI's, yet in such
 cases it loses no more than $0.68\%$.
Over most of the parameter region the MPI policy is
 better, achieving maximum performance gains of over
$7.6\%$.

\subsection{Scheduling a multiclass make-to-order/make-to-stock
     queue}
   \label{s:smmtomtsq}
A fundamental problem arising in manufacturing systems is that of 
dynamic scheduling of production in a multiproduct system, which is
conveniently modeled through a multiclass make-to-stock  queue that
incorporates stock and backorder holding costs.

\citet{wein92} addresses the problem of designing a dynamic scheduling policy
to minimize average cost, by carrying out a heavy-traffic analysis
based on solving an approximating Brownian control  problem.
The proposed policy emerging from such an analysis has the following
form. 
The server or machine is idle as long as the weighted inventory
process defined in that paper is above a certain threshold level, and
no classes are in ``danger of being backordered,'' as explained there.
Otherwise, the machine is allocated to some class, according to the
following static index rule. 
If there is some class in danger of being backordered, 
the machine is assigned to a class of \emph{largest index} $b_k \mu_k$, where $b_k$
and $\mu_k$ are the backorder cost rate and the service rate for class
$k$, respectively.
If such is not the case, the machine is assigned to a class of
\emph{smallest index} $h_k \mu_k$, where $h_k$ is class $k$'s stock
holding cost rate.

\citet{vewe} investigate several policies for such a
model, based on identifying a hedging point of base stock levels,
which determines the idleness region in state space,  along
with an index policy that determines the machine allocation. 
Though they consider use of the restless bandit index policy in the
lost-sales case
where no backorders are allowed, they discard such an approach in
the case with backorders, as the index does not exist under the 
average criterion (cf.\ Sect.\ \ref{s:aoc} above). 
They thus propose a policy based on a heavy-traffic diffusion approximation.

\citet{pezipk97} propose several dynamic index policies based on intuitive
limited look-ahead arguments. 
They report experimental results showing that the policy  they
term \emph{myopic$(T)$}, which uses as the look-ahead time
determining the index of a class the
sojourn time of a job, appears both to be nearly optimal and to
outperform alternative index policies.
\citet{vekada00} furnish theoretical support for use of such an
index, showing that it partially characterizes optimal policies in a
Markovian model.

We must also mention the work in 
\citet{dushongl}, which applies restless bandit indexation to address
a Brownian model of a multiclass make-to-stock system.

\citet{nmmor06} develops results announced in \citet{nmaller03} on
 indexation analysis for $M/G/1$ make-to-order/make-to-stock queues,
 aimed at deriving index policies for  a multiclass $M/G/1$ queueing model where some classes
 must be processed in make-to-order mode, whereas others may be
 processed in
 make-to-stock mode. The model incorporates state-dependent stock
 holding and backorder cost rates for each class.
For earlier work and applications of such combined system models see,
 e.g.,  \citet{adanwal98}, \citet{somanetal}, and the references
 therein.
Assuming convexity of backorder and holding cost rates, it is shown in
 \citet[Sect.\ 6]{nmmor06} that the restless bandit of concern, which are single
  make-to-order/make-to-stock queues subject to service control, are
 PCL-indexable and hence indexable, both under the discounted
 criterion and under the mixed average-bias criterion reviewed in
 Sect.\ \ref{s:aoc}.

Taking as the state of a make-to-order/make-to-stock  $M/G/1$ queue the net backorder level 
$X(t) = X^+(t) - X^-(t)$, where $X^+(t)$ is the number of backordered items
and $X^-(t)$ is the number of finished items in stock at each time
$t$, we may conveniently represent the backorder and stock holding
cost rate by a single state-dependent cost rate $h_i$.
The average-bias MPI for such a system is then given by the simple
expression (where $\Delta h_i \triangleq h_i - h_{i-1}$)
\begin{equation}
\label{eq:mtsi}
\nu_i^* = \mu \Ex\left[\Delta h_{X + i}\right],
\end{equation}
where $\Delta h_i \triangleq h_i - h_{i-1}$, $\mu$ is the service rate and $X$ is a random variable having the steady-state distribution of
the number-in-system for a corresponding make-to-order $M/G/1$ queue.

Consider the special case where backorder and stock holding
cost rates are linear, so that 
\[
h_i = 
\begin{cases}
b i & \text{if } i \geq 1 \\
- h i & \text{otherwise},
\end{cases}
\]
where $b$ (resp.\ $h$) is the backorder (resp.\ stock) holding cost
rate per item per unit time.
Then, we show in eq.\ (42) of the paper referred to above that the MPI has the
simple evaluation
\begin{equation}
\label{eq:linmpimtso}
\nu_i^* = 
\begin{cases}
b \mu & \text{if } i \geq 1 \\
\left[\left(b + h\right) 
  \Prob\left\{X \geq -i+1\right\} - h\right] \mu & 
\text{otherwise}.
\end{cases}
\end{equation}

Now, notice that the Pe\~na-P\'erez and Zipkin myopic$(T)$ index is
formulated in the present notation by
\[
\nu^{\textup{PPZ}}_i = 
\begin{cases}
b \mu & \text{if } i \geq 1 \\
\left[\left(b + h\right)
  \Prob\left\{D(T) \geq -i+1\right\} - h\right] \mu, & 
\text{otherwise}.
\end{cases}
\]
where $D(T)$ has the steady-state distribution of the number of
arrivals during a customer's sojourn
time in a corresponding make-to-order $M/G/1$ queue under FIFO.

It turns that the indices $\nu^{\textup{PPZ}}_i$ and $\nu_i^*$ above
are identical. The reason is that
$X$ and $D(T)$ have the same distribution, as argued
in \citet[Chapt.\ 5.7]{kleinr75}. Notice that such a result is a form of the
distributional form of Little's law. See \citet{hajinew71} and \citet{keilserv88}.
Thus, the MPI recovers the Pe\~na-P\'erez and Zipkin myopic$(T)$ index
in the case of linear costs, and further extends it to models with convex nonlinear cost
rates and also to the discounted criterion.

Denoting the $\nu_i^{*, \alpha}$ the discounted MPI in the $M/M/1$
case with discount rate
$\alpha > 0$, it is  further of interest to evaluate the following limiting 
\emph{myopic index}:
\[
\nu_i^{\mathrm{myopic}} \triangleq \lim_{\alpha \nearrow \infty} 
\alpha \nu_i^{*, \alpha} = \mu \Delta h_i = \begin{cases}
b \mu & \text{if } i \geq 1 \\
 - h \mu & 
\text{otherwise}.
\end{cases}
\]
As pointed out at the end of \citet{nmmor06}, such a myopic
index is precisely the index obtained by \citet{wein92} through a
heavy-traffic diffusion analysis.

\subsection{Scheduling a multiclass queue with finite buffers}
\label{s:smqfb}
A remarkable paradox in queueing theory is that, 
while in all queueing systems arising in the real world the
storage or buffer space for holding customers is limited, most queueing models 
assume the latter to be unbounded.
Such a state of affairs is probably related to the fact that
 infinite-buffer models are typically more tractable than their
 finite-buffer counterparts, as the latter exhibit complex boundary
 effects that complicate the analyses. 
The research on models for scheduling multiclass queues is no
exception. Thus, e.g.,  the optimality of the classic $c \mu$ and the Klimov index
rules is established assuming unlimited storage
capacity.

In contrast, limited attention has been devoted to problems
involving the dynamic scheduling of
multiclass queues where each class has its own finite dedicated buffer
space, despite their relevance in manufacturing (scheduling of parts
for processing at a flexible machine) and computer-communication (scheduling of
packets for access to a transmission channel) systems.
Finite buffers introduce the possibility of blocking, which occurs when
customers arrive to find a full buffer and are hence rejected and lost.
In such systems, optimal policies have only been identified under
strong symmetry conditions. 
Thus, \citet{spar93} addressed the  dynamic scheduling
problem with the goal of minimizing the expected blocking cost, in a
Markovian model where all classes have the same arrival, service and blocking
cost rates, whereas buffer sizes may differ. 
They exploited a duality relationship between routing and scheduling
problems with finite buffers to infer, from the known optimality of
the \emph{shortest nonfull queue} routing rule established in
\citet{hk90}, the optimality of its dual scheduling rule, which
prescribes to dynamically allocate the server at each time to a nonempty class with
the \emph{smallest residual capacity} (SRC). 
Such a result was extended to the case of general service-time 
distributions in \citet{wassbamb96}.

\citet{kimvo98} addressed a nonsymmetric Markovian model having two classes $k$
with holding and rejection cost rates $c_k$ and $r_k$, respectively.
They identified a condition under which discount-optimal policies are
characterized by a monotonic switching curve so that, if it is
optimal to serve class $1$ in state $(i_1, i_2)$, then it is also
optimal to serve it in state $(i_1+1, i_2)$, and correspondingly for
class $2$, where the state is the number of customers in each class.
The required condition is that, for each class $k$, 
\begin{equation}
\label{eq:ckvo}
\alpha r_k \geq c_k,
\end{equation}
where $\alpha > 0$ is the discount rate.
Notice that (\ref{eq:ckvo}) means that the cost $r_k$ of rejecting a
customer is greater than or equal to the cost $c_k / \alpha$ of
holding it forever in the system.
They further show by example that, if such a condition is violated,
optimal policies can exhibit complex patterns, pointing out the
seemingly counterintuitive phenomenon that, in
certain situations, there may be an incentive to award higher priority
to a shorter queue.

\citet{nmqs06} deploys restless bandit indexation theory to a Markovian
multiclass queue with dedicated finite buffers, finding that condition
(\ref{eq:ckvo}) plays a key role in the analyses and results.
It turns out that the state ordering induced by the MPI in a traffic class
depends critically on whether or not such a condition holds.
If it does, let us say that the class is
\emph{loss-sensitive}. Otherwise, we say that it is \emph{delay-sensitive},
in each case relative to the prevailing discount rate $\alpha$.
Such concepts are particularly relevant in the setting of contemporary 
research efforts to provide differentiated service to
 heterogeneous traffic streams in
packet-switched computer-communication networks, where some traffic 
classes are more sensitive to losses (e.g., file transfer), while 
others are more sensitive to delays (e.g., interactive and multimedia 
applications).

That paper shows that a loss-sensitive class is PCL-indexable,
relative to the nested active-set family (cf.\ Sect.\ \ref{s:awrv}) 
consistent with the index ordering
\begin{equation}
\label{eq:ordls}
\nu_n^{\alpha, *} \leq \nu_{n-1}^{\alpha, *} \leq \cdots \leq \nu_1^{\alpha, *},
\end{equation}
where the state $i$ in $\nu_i^{\alpha, *}$ is now taken to be the
\emph{number of empty buffer spaces}, as it does not depend on the
buffer size $n$,  and the notation makes explicit the
index' dependence on the discount rate. Notice that the ordering in (\ref{eq:ordls})
agrees with the intuition that, other things being equal,
queues with fewer empty buffer spaces should be awarded higher service priority.

In the pure loss-sensitive case $r > 0 = c$ (where we have dropped the
class label $k$ from the notation), 
the natural index candidate for  scheduling
under the average criterion is obtained as the
limiting index obtained from $\nu_i^{\alpha, *}$ as $\alpha$
vanishes.
The latter turns out to be constant, being equal to $r \mu$, which is
noninformative in the case where all rejection and service rates are
the same. To break ties, we consider the Maclaurin series expansion 
$\nu_i^{\alpha, *} = r \mu - \alpha \gamma_i^* + o(\alpha)$ as $\alpha
\searrow 0$, which yields the \emph{second-order MPI}
\begin{equation}
\label{eq:gammak}
\gamma_i^* = 
\begin{cases}
\displaystyle{\frac{r}{\rho} 
\left\{i + 1 + \frac{1/\rho^{i} - (1-\rho) i - 1}{(1-\rho)^2}\right\}} & \text{if } \rho \neq 1 \\ \\
\displaystyle{r \frac{(i+1) (i + 2)}{2}} & \text{if } \rho = 1,
\end{cases}
\end{equation}
where $\rho \triangleq \lambda/\mu$.
Thus, among classes with the same first-order MPI $r \mu$, 
the resultant index rule gives higher service priority to classes with \emph{smaller} values of the
second-order MPI $\gamma_i^*$.
We must remark that, in the classic Bayesian multiarmed
bandit problem with Bernoulli arms,  a tie-breaking second order index was
shown in \cite{kellymab} to give an optimal policy for values of the discrete-time
discount factor near one.

The paper further shows that a delay-sensitive class is PCL-indexable,
relative to the nested active-set family consistent with the ordering
(\ref{eq:ordls}), which is however interpreted by taking now the state
$i$ in $\nu_i^{\alpha, *}$ to represent the
\emph{number of customers in the system}.
Hence, other things being equal, such an ordering prescribes to award
higher priority to \emph{shorter nonempty queues}.
While such a result might appear at first sight as counterintuitive,
it agrees with the observation of Kim and Van Oyen referred to above,
and also with experimental results on the structure of optimal
policies for particular instances.
One possible interpretation is that, since holding
costs are bounded, when such costs are dominant it is more productive
to prevent congestion than to react to it, while the opposite holds
otherwise.

Note further that a class with positive holding cost rate $c > 0$ is
always delay-sensitive for small enough values of the discount rate
$\alpha$, which leads us to consider the limiting index obtained
as the latter vanishes. This has the evaluation
\begin{equation}
\label{eq:lampi}
\nu_i^{*} =
\begin{cases}
\displaystyle \frac{c}{\rho} \left\{n - \frac{\rho}{1-\rho} + i
  \frac{\rho^{i}}{1-\rho^{i}}\right\}   + 
 r \mu & \text{ if } \rho \neq 1
  \\ \\
\displaystyle c \left\{n-\frac{i-1}{2}\right\} + r \mu & 
\text{ if } \rho = 1.
\end{cases}
\end{equation}
Such a limiting index is shown to be indeed an MPI relative to the bias 
criterion discussed in that paper, which is appropriate as a
priority-index for the corresponding multiclass system. 

Yet, Are such MPI policies useful?
We refer the reader to the experimental results on two-class instances
reported in
\citet{nmqs06}. 
Across the instances investigated, such policies are typically near
optimal, and outperform, often substantially, traditional naive
policies such as the $c \mu$ rule or the SRC rule referred to above.  
We only identified a range of instances where the MPI policies appear to perform poorly:
in the pure loss-sensitive case with distinct first-order indices 
$r_k \mu_k$, which yields a static priority policy.

\section{More recent work}
\label{s:rw}
\subsection{Algorithmic characterization of indexability}
\label{s:aci}
\citet{nmbfsind07} announces results of deploying a parametric LP approach, based on the
classic parametric-objective simplex algorithm of \citet{gassaaty55},
to test for indexability and compute the MPI of an indexable
semi-Markovian restless
bandit instance. 
The resultant \emph{Complete-Pivoting Indexability} (CPI) algorithm
consists of an initialization stage, which computes the initial simplex
tableau, followed by a loop that 
performs $2 n^3 + O(n^2)$ arithmetic operations for an indexable
project with $n$ controllable states.
The algorithm is given in a detailed
\emph{block-partitioned} form, i.e., based on operations on
submatrices (blocks) of a base matrix, which is ready for actual
implementation. 
The importance of such block implementation has been emphasized in the
scientific computing literature, where it is advocated as a means to
partially overcome the exponentially widening gap between memory access
times and processor speed in contemporary computers. 
See, e.g., \citet{doneij00}.

We also present in that paper a \emph{Reduced-Pivoting Indexability}
(RPI) algorithm, that avoids some unnecessary operations in pivot
steps to achieve a
reduced operation count of $n^3 + O(n^2)$ operations in its main loop.
Yet, such an improved theoretical complexity is obtained at the
expense of manipulating submatrices with arbitrary row and column indices,
which results in relatively inefficient strided memory-access
patterns. 
The latter are known to produced significant slowdowns in actual
running times, which is verified experimentally in the paper. Despite
the above operation counts, the CPI algorithm is substantially faster
in practice than the RPI algorithm. 

We have used such algorithms, in MATLAB implementations developed by
the author,  to assess empirically via simulation the prevalence of 
indexability and PCL-indexability on randomly generated restless
bandit instances with dense transition probability matrices. 
Specifically, in each instance active rewards and transition
probabilities were generated with MATLAB as $\textup{Uniform}(0, 1)$
pseudo-random numbers, and then transition matrices were properly
scaled dividing each row by its sum. 
For each of the state space sizes $n = 3, \ldots, 7$ a sample of 
$10^7$ instances was drawn and, for each instance, the discount factor 
$\beta$ was varied from $0.1$ to $1$. Table \ref{tab:crsmpi} reports
the counts obtained of nonindexable instances and of indexable yet
non-PCL instances for each $(\beta, n)$ combination.
The results suggest that both indexability and PCL-indexability are
highly prevalent properties, with their incidence sharply increasing as the 
discount factor gets smaller and as the state space gets larger.

\begin{table}
\caption{Counts on samples of $10^7$ restless bandit instances.}
\begin{center}
\begin{tabular}{c|rrrrr|rrrrr} 
& \multicolumn{5}{c|}{Nonindexable} & \multicolumn{5}{c}{Indexable
  non-PCL} \\ \hline
& \multicolumn{5}{c|}{number of states} & \multicolumn{5}{c}{number of
 states} \\
$\beta$ & $3$ & $4$ & $5$ & $6$ & $7$  & $3$ & $4$ & $5$ & $6$ & $7$ \\
\hline
$0.1$   & $0$ & $0$ & $0$ & $0$ & $0$  & $0$ & $0$ & $0$ & $0$ & $0$
\\
$0.2$   & $0$ & $0$ & $0$ & $0$ & $0$  & $0$ & $0$ & $0$ & $0$ & $0$
\\
$0.3$   & $0$ & $0$ & $0$ & $0$ & $0$  & $0$ & $0$ & $0$ & $0$ & $0$
\\
$0.4$   & $0$ & $0$ & $0$ & $0$ & $0$  & $0$ & $0$ & $0$ & $0$ & $0$
\\
$0.5$   & $0$ & $0$ & $0$ & $0$ & $0$  & $0$ & $0$ & $0$ & $0$ & $0$
\\
$0.6$   & $0$ & $0$ & $0$ & $0$ & $0$  & $0$ & $0$ & $0$ & $0$ & $0$
\\
$0.7$   & $0$ & $0$ & $0$ & $0$ & $0$  & $30$ & $0$ & $0$ & $0$ & $0$
\\
$0.8$   & $16$ & $1$ & $0$ & $0$ & $0$  & $574$ & $32$ & $1$ & $0$ & $0$
\\
$0.9$   & $135$ & $7$ & $0$ & $0$ & $0$  & $4460$ & $509$ & $36$ & $5$ & $0$
\\
$1.0$   & $818$ & $66$ & $4$ & $0$ & $0$  & $18631$ & $3640$ & $425$ & $50$ & $3$ 
\end{tabular}
\end{center}
\label{tab:crsmpi}
\end{table}

\begin{figure}
\centering
\begin{psfrags}%
\psfragscanon%
\psfrag{s03}[l][l]{\color[rgb]{0,0,0}\setlength{\tabcolsep}{0pt}\begin{tabular}{l}$\emptyset$\end{tabular}}%
\psfrag{s04}[l][l]{\color[rgb]{0,0,0}\setlength{\tabcolsep}{0pt}\begin{tabular}{l}$\{1\}$\end{tabular}}%
\psfrag{s05}[c][l]{\color[rgb]{0,0,0}\setlength{\tabcolsep}{0pt}\begin{tabular}{l}$\{2\}$\end{tabular}}%
\psfrag{s06}[l][l]{\color[rgb]{0,0,0}\setlength{\tabcolsep}{0pt}\begin{tabular}{l}$\{3\}$\end{tabular}}%
\psfrag{s07}[l][l]{\color[rgb]{0,0,0}\setlength{\tabcolsep}{0pt}\begin{tabular}{l}$\{1, 2\}$\end{tabular}}%
\psfrag{s08}[l][l]{\color[rgb]{0,0,0}\setlength{\tabcolsep}{0pt}\begin{tabular}{l}$\{1, 3\}$\end{tabular}}%
\psfrag{s09}[c][l]{\color[rgb]{0,0,0}\setlength{\tabcolsep}{0pt}\begin{tabular}{l}$\{2, 3\}$\end{tabular}}%
\psfrag{s10}[l][l]{\color[rgb]{0,0,0}\setlength{\tabcolsep}{0pt}\begin{tabular}{l}$\{1, 2, 3\}$\end{tabular}}%
\psfrag{s11}[t][t]{\color[rgb]{0,0,0}\setlength{\tabcolsep}{0pt}\begin{tabular}{c}$g^{\pi}$\end{tabular}}%
\psfrag{s12}[b][b]{\color[rgb]{0,0,0}\setlength{\tabcolsep}{0pt}\begin{tabular}{c}$f^{\pi}$\end{tabular}}%
%
\psfrag{x01}[t][t]{0}%
\psfrag{x02}[t][t]{0.1}%
\psfrag{x03}[t][t]{0.2}%
\psfrag{x04}[t][t]{0.3}%
\psfrag{x05}[t][t]{0.4}%
\psfrag{x06}[t][t]{0.5}%
\psfrag{x07}[t][t]{0.6}%
\psfrag{x08}[t][t]{0.7}%
\psfrag{x09}[t][t]{0.8}%
\psfrag{x10}[t][t]{0.9}%
\psfrag{x11}[t][t]{1}%
%
\psfrag{v01}[r][r]{0}%
\psfrag{v02}[r][r]{0.1}%
\psfrag{v03}[r][r]{0.2}%
\psfrag{v04}[r][r]{0.3}%
\psfrag{v05}[r][r]{0.4}%
\psfrag{v06}[r][r]{0.5}%
\psfrag{v07}[r][r]{0.6}%
\psfrag{v08}[r][r]{0.7}%
\psfrag{v09}[r][r]{0.8}%
\psfrag{v10}[r][r]{0.9}%
\psfrag{v11}[r][r]{1}%
%
\includegraphics[height=2.2in,width=6.1in,keepaspectratio]{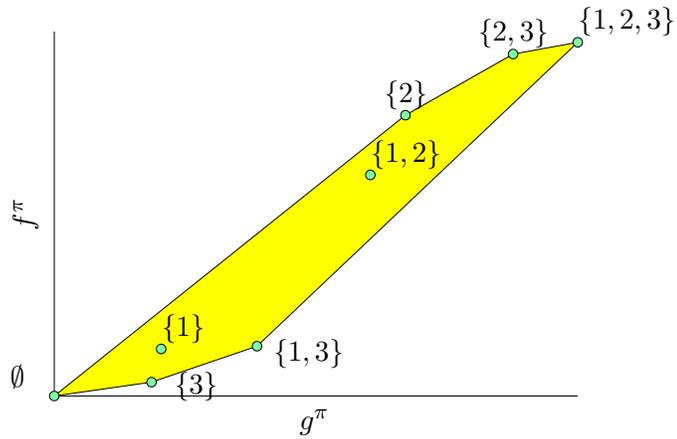}%
\end{psfrags}%
\caption{Achievable work-reward performance region of an indexable non-PCL bandit.}
\label{fig:indnpcl3}
\end{figure}

To illustrate by a concrete example how a bandit can be indexable yet
 not PCL-indexable, consider Figure \ref{fig:indnpcl3}, which displays
 the achievable work-reward performance
 region (cf.\ Sect.\ \ref{s:awrv}) for the three-state instance with $\beta = 0.9$, 
\[
\mathbf{R}^1 =
\begin{bmatrix}
0.44138 \\ 0.8033 \\ 0.14257
\end{bmatrix}, 
\mathbf{P}^1 =
\begin{bmatrix}
0.1719 & 0.1749 & 0.6532 \\
0.0547 & 0.9317 & 0.0136 \\
0.1547 & 0.6271 & 0.2182
\end{bmatrix}, 
\mathbf{P}^0 =
\begin{bmatrix}
0.3629 & 0.5028 & 0.1343 \\
0.0823 & 0.7534 & 0.1643 \\
0.2460 & 0.0294 & 0.7246
\end{bmatrix},
\]
and $\mathbf{R}^0 = \mathbf{0}$.
The plot shows that this is an indexable instance, relative to the nested
family of optimal policies $\mathcal{F}_0 = \{\emptyset, \{2\},
\{2, 3\}, \{1, 2, 3\}\}$.
Yet, it is not PCL-indexable, since $g^{\{1, 2\}} < g^{\{2\}}$ and
hence (cf.\ (\ref{eq:equivrefwisp})), $w_1^{\{2\}} < 0$.
In contrast, Figure \ref{fig:indpcl0} shows the achievable work-reward
performance region of a PCL-indexable project.

\subsection{More powerful indexability conditions and faster index
     computation}
\label{s:mpicfic}
When faced with a particular restless bandit model, it is of interest
for researchers to have practical methods to establish analytically its
indexability under an appropriate parameter range. 
While we developed the PCL-indexability approach discussed
in Sect.\ \ref{s:sicaga} for such a purpose, our more recent work has
revealed limitations that have prompted further research. 
In particular: (i) we have found that analytical verification of the
second condition in 
Definition \ref{def:pclib}, namely that the sequence of index
values produced by the adaptive-greedy algorithm is monotone
nonincreasing, can be overly hard or even elusive in models with a
multi-dimensional state, such as that discussed in Sect.\ \ref{s:mabd};
and (ii) we have further found, as discussed in Sect.\ \ref{s:mabsd}, a
relevant restless bandit model that is indexable, yet
not PCL-indexable, as the first condition in Definition
\ref{def:pclib} does not necessarily hold. 

Such a state of affairs motivated the author to introduce in 
\citet{nmbfsind07} sufficient
indexability conditions that are both easier to apply and less
restrictive than PCL-indexability. 
We next review such conditions, termed \emph{LP-indexability}
conditions as they
are based on LP analyses.
As the PCL-indexability conditions in Sect.\ \ref{s:sicaga}, the new
conditions are also based on positing an appropriate set system
$(N^{\{0, 1\}}, \mathcal{F})$. 
Yet, stronger requirements are imposed on the latter than those in 
Assumption \ref{ass:nfreq}, based on the concept of 
\emph{monotonically connected set
  system} introduced in that paper. 

\begin{assumption}
\label{def:ass2} 
$(N^{\{0, 1\}}, \mathcal{F})$ is  a \emph{monotonically connected set
  system}, i.e., it satisfies:
\begin{itemize}
\item[\textup{(i)}] $\emptyset, N^{\{0, 1\}} \in \mathcal{F}$;
\item[\textup{(ii)}] for every 
$S, S' \in \mathcal{F}$ with $S \subset S'$ there exist 
$j \in \partial^{\textup{out}}_{\mathcal{F}} S$ and $j' \in
\partial^{\textup{in}}_{\mathcal{F}} S'$ such that $S \cup
\{j\} \subseteq S'$ and $S \subseteq S' \setminus \{j'\}$; and
\item[\textup{(iii)}] for any $S, S' \in \mathcal{F}$,
  $S \cup S' \in \mathcal{F}$ and $S \cap S' \in \mathcal{F}$.
\end{itemize}
\end{assumption}

The term ``monotonically connected'' is motivated by the
fact that, in such a set system, one can always
connect two feasible sets $S \subset S'$ by a monotone increasing
sequence $S_1 \subset \cdots \subset S_m$ of adjacent sets in $\mathcal{F}$,
with $S_1 = S$, $S_m =
S'$. Further, one can also connect two distinct feasible sets $S
\neq S'$ through two successive monotone sequences of adjacent sets in
$\mathcal{F}$, the first of which is monotone increasing and connects
$S$ to $S \cup S'$, while the
second is monotone decreasing and connects $S \cup S'$ to $S'$.

We further write
\begin{equation}
\label{eq:lud}
\underline{r}^S \triangleq \max_{j \in S^c, w_j^S = 0} r_j^S \quad
\text{ and } \quad 
\overline{r}^S \triangleq \min_{j \in S, w_j^S = 0} r_j^S,
\end{equation}
 adopting the convention that the maximum (resp. minimum) over
an empty set has the value $-\infty$ (resp. $+\infty$).

\begin{definition}[LP$(\mathcal{F})$-indexability]
\label{def:bfsi} 
We say that a restless bandit is \emph{LP$(\mathcal{F})$-indexable}
  if:
\begin{itemize}
\item[(i)] $w_i^\emptyset, w_i^{N^{\{0, 1\}}} \geq 0$ for $i \in
 N^{\{0, 1\}}$,
 and $\underline{r}^\emptyset  \leq 0 \leq \overline{r}^{N^{\{0, 1\}}}$;
\item[(ii)]
for each $S \in \mathcal{F}$,  
$w_i^S > 0$ for $i \in
   \partial^{\textup{in}}_{\mathcal{F}} S
   \cup \partial^{\textup{out}}_{\mathcal{F}} S$; and
\item[(iii)] for every wage $\nu \in \mathbb{R}$ there exists an optimal
  active set $S \in \mathcal{F}$ for
  (\ref{eq:nuwp}).
\end{itemize}
\end{definition}

In practice, to establish that a bandit model is
LP$(\mathcal{F})$-indexable one would prove conditions
(i) and (ii) in Definition \ref{def:bfsi} by carrying out an 
analysis of relevant marginal work and reward measures.
As for condition (iii), it requires us to prove a structural property
of optimal policies for $\nu$-wage problem (\ref{eq:nuwp}), for which 
DP techniques can often be useful. 

The interest of such a bandit class is due to the following result. 
\begin{theorem}
\label{the:bfsfi} The following holds\textup{:} 
\begin{itemize}
\item[\textup{(a)}]
An LP$(\mathcal{F})$-indexable bandit is indexable, and its MPI is
computed in nonincreasing order by adaptive-greedy algorithm
$\mathrm{AG}_{\mathcal{F}}$ in Table \textup{\ref{fig:agaf}}.
\item[\textup{(b)}] An indexable bandit is
  LP$(\mathcal{F}_0)$-indexable relative to some nested active-set family 
 $\mathcal{F}_0$.
\end{itemize}
\end{theorem}

Note that part (a) of Theorem \ref{the:bfsfi} yields a sufficient condition for
indexability, and further extends the scope of the adaptive-greedy
algorithm in Table \textup{\ref{fig:agaf}} from the class of
PCL-indexable bandits to the wider class of LP-indexable bandits,
i.e., those that are LP$(\mathcal{F})$-indexable relative to some 
family $\mathcal{F}$. 
Part (b) shows that the latter class encompasses all indexable
bandits.

For a specific application of Theorem \ref{the:bfsfi} we refer the
reader to \citet{nmswd07}, which analyzes a bandit model that is
LP-indexable yet not PCL-indexable. 
That paper further shows that the condition (ii) in Definition \ref{def:pclib}
on PCL$(\mathcal{F})$-indexability can be replaced by condition (iii) in 
Definition (\ref{def:bfsi}).

The validity of index algorithm 
$\mathrm{AG}_{\mathcal{F}}$ for LP$(\mathcal{F})$-indexable bandits
allows us to leverage structural knowledge on optimal policies to
compute the index faster. 
Thus, a fast-pivoting block implementation is given in \citet{nmbfsind07} of
such an algorithm which, after an initialization stage involving the
solution of a block linear equation system, performs
in its main loop $(2/3) n^3 + O(n^2)$ arithmetic operations for a
project having
$n$ controllable states. 
Such an algorithm extends that introduced in \citet{nmijoc06} for
computing the Gittins index and solving the problem of optimal
stopping of a finite Markov chain, which has the same complexity in its
loop yet does not require the initialization stage. 
Notice that $(2/3) n^3 + O(n^2)$ is the complexity of solving an $n \times n$
linear equation system by Gaussian elimination. 

\begin{figure}[htb]
\centering
\begin{psfrags}%
\psfragscanon%
%
\psfrag{s01}[b][b]{\color[rgb]{0,0,0}\setlength{\tabcolsep}{0pt}\begin{tabular}{c}\end{tabular}}%
\psfrag{s02}[t][t]{\color[rgb]{0,0,0}\setlength{\tabcolsep}{0pt}\begin{tabular}{c}$n$\end{tabular}}%
\psfrag{s03}[b][b]{\color[rgb]{0,0,0}\setlength{\tabcolsep}{0pt}\begin{tabular}{c}runtime (hours)\end{tabular}}%
%
\psfrag{x01}[t][t]{0}%
\psfrag{x02}[t][t]{0.1}%
\psfrag{x03}[t][t]{0.2}%
\psfrag{x04}[t][t]{0.3}%
\psfrag{x05}[t][t]{0.4}%
\psfrag{x06}[t][t]{0.5}%
\psfrag{x07}[t][t]{0.6}%
\psfrag{x08}[t][t]{0.7}%
\psfrag{x09}[t][t]{0.8}%
\psfrag{x10}[t][t]{0.9}%
\psfrag{x11}[t][t]{1}%
\psfrag{x12}[t][t]{1000}%
\psfrag{x13}[t][t]{2000}%
\psfrag{x14}[t][t]{3000}%
\psfrag{x15}[t][t]{4000}%
\psfrag{x16}[t][t]{5000}%
\psfrag{x17}[t][t]{6000}%
\psfrag{x18}[t][t]{7000}%
\psfrag{x19}[t][t]{8000}%
\psfrag{x20}[t][t]{9000}%
\psfrag{x21}[t][t]{10000}%
%
\psfrag{v01}[r][r]{0}%
\psfrag{v02}[r][r]{0.1}%
\psfrag{v03}[r][r]{0.2}%
\psfrag{v04}[r][r]{0.3}%
\psfrag{v05}[r][r]{0.4}%
\psfrag{v06}[r][r]{0.5}%
\psfrag{v07}[r][r]{0.6}%
\psfrag{v08}[r][r]{0.7}%
\psfrag{v09}[r][r]{0.8}%
\psfrag{v10}[r][r]{0.9}%
\psfrag{v11}[r][r]{1}%
\psfrag{v12}[r][r]{0.5}%
\psfrag{v13}[r][r]{1}%
\psfrag{v14}[r][r]{1.5}%
\psfrag{v15}[r][r]{2}%
\psfrag{v16}[r][r]{2.5}%
\psfrag{v17}[r][r]{3}%
\psfrag{v18}[r][r]{3.5}%
\psfrag{v19}[r][r]{4}%
\psfrag{v20}[r][r]{4.5}%
%
\includegraphics[height=2.8in,width=6.1in,keepaspectratio]{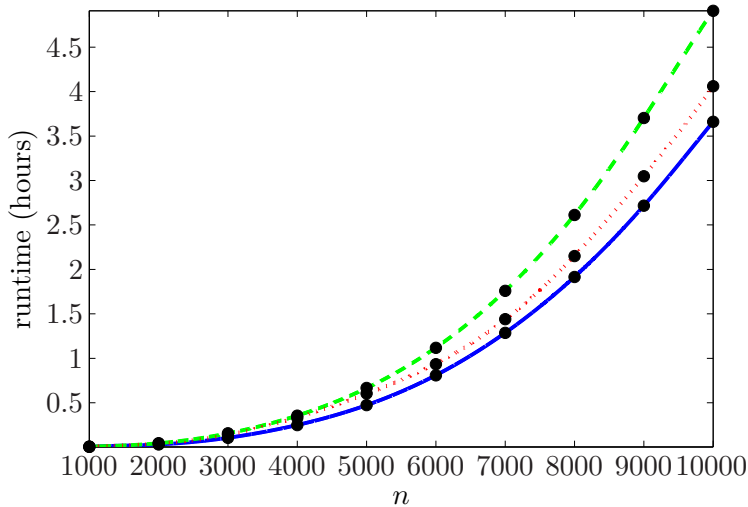}%
\end{psfrags}%
\caption{Runtimes with cubic least-squares fit: FPAG (solid), CPI
  (dotted) and RPI (dashed) algorithms.}
\label{fig:rclsf}
\end{figure}

To give an idea of actual runtimes of the three algorithms discussed, we report the results of an
experiment based on the author's MATLAB implementations.
A random restless bandit instance was generated for each of the state
space sizes $n = 1000$ up to $10000$ in increments of $1000$, with a
discount factor $\beta = 0.8$.
Each instance was tested both for indexability and PCL-indexability,
yielding a positive result in each case. 
Figure \ref{fig:rclsf} shows the runtime performance of each
algorithm, where FPAG refers to the fast-pivoting implementation of the
adaptive-greedy algorithm mentioned above.
The experiment was conducted on an HP xw9300 dual-processor Opteron
254 (2.8 GHz) workstation running MATLAB 2007a under Windows XP x64.
The new hyperthreading capability of MATLAB in that release was
enabled to allow computations to use both system processors.

The results show that the FPAG algorithm is the faster of the three,
followed by the CPI and then the RPI algorithms. 
The relative runtime performance does not match what would be expected
based on theoretical operation counts. 
Such a phenomenon is well known in the field of scientific computing,
as memory-access patterns are often the dominant factor in an
algorithm's actual performance. 
See, e.g., \citet{doneij00}. 
Thus, algorithms RPI and FPAG achieve their
reduced operation counts at the expense of working on submatrices with 
complex row and index patterns, which results in strided access to
memory, whereas algorithm CPI's work is concentrated on a whole matrix
manipulated as a single contiguous block. 

\subsection{Scheduling a multiclass wireless queue with finite
     buffers}
\label{s:smwqfb}
The ubiquitous presence of wireless packet-switched 
computer-communication networks
motivates the investigation of suitable modifications of traditional
dynamic scheduling models.
Thus, e.g., in a system scheduling multiple traffic streams
generated by mobile users that vie for 
access to a transmission channel, 
the wireless environment creates the possibility that the latter is
temporarily unavailable to some users, due to phenomena such as
fading. 
A convenient model for such a system is given by a multiclass queue with
a server whose connectivity to each customer class is turned on an 
off in a random fashion, according to a two-state Markov chain.

Variants of such a model have attracted considerable research
interest, mostly aimed at establishing optimality of greedy scheduling
policies under rather stringent conditions on model parameters.
See, e.g., 
\citet{taseph93}, \citet{lotttenek00}, and \citet{bambmich02}.

\citet{nmngi06} announces results of deploying
 restless bandit indexation to obtain 
new dynamic scheduling policies in a
discrete-time Markovian 
multiclass wireless queue with finite dedicated buffers.
Specifically, we assume that there are 
$K$ distinct traffic classes,
 labeled by $k \in \mathbb{K} \triangleq \{1, \ldots, K\}$. 
During a time period (or \emph{slot}
 in computer-communications lingo) the number of 
class $k$ customers that arrive is distributed as a Poisson random
variable with rate $\lambda_k$.
If the server is allocated to a class $k$ customer during a period, 
the probability that the service is completed by the end of the period
is $\mu_k$.
We thus assume Poisson arrival 
processes and geometric service times, which are
mutually independent across periods and classes.
Upon a class $k$ customer's arrival, it joins its class' queue,
which is held in a dedicated finite buffer of size $n_k$, if this is not full;
otherwise, the customer is blocked and lost.
We denote by $X_k(t)$ the number of class $k$
customer in system at the start of period $t$. 

The time-varying connectivity 
is modeled through a binary \emph{connectivity process} $s_k(t)$ for each 
class $k$, where $s_k(t) = 1$ when class $k$ is connected to the 
server during period $t$, and $s_k(t) = 0$ otherwise.
We assume that each $s_k(t)$ evolves as a binary ergodic Markov chain with
transition probabilities $p_k(e_1, e_2)$, for $e_1, e_2 \in \{0, 1\}$, 
and that connectivity processes are independent 
across classes.

We take the \emph{state} of class $k$ at time $t$ as 
the connectivity-number of customers pair $(s_k(t), X_k(t))$.
Its \emph{state space} is thus 
\[
N_k \triangleq \{0, 1\} \times \{0, \ldots, n_k\}.
\]

A controller chooses at the start of each period
the nonempty connected class, if any, to be served. 
Customers within a class are serviced in FIFO order.
Such choices are represented by binary \emph{action processes}
$a_k(t)$, where $a_k(t) = 1$ if the server is working on class $k$ in
period $t$ and $a_k(t) = 0$ otherwise. The corresponding sample-path 
service capacity constraint is thus
\[
\sum_{k \in \mathbb{K}} a_k(t) \leq 1, \quad t \geq 0.
\]
Notice that, in states $(s_k(t), X_k(t))$ with $s_k(t) = 0$ or $X_k(t) = 0$ the only meaningful action is 
$a_k(t) = 0$. 
We will thus consider such states, whose union we denote by
$N_k^{\{0\}}$, as \emph{uncontrollable}. 
The remaining set of class $k$ states, which we denote by $N_k^{\{0, 1\}}$, are
\emph{controllable}, in that both actions are available in them and differ.

Action choice is dynamically 
prescribed through adoption of a 
\emph{scheduling policy} $\pi$. 
This is 
 chosen from the space $\Pi$ of \emph{admissible policies}, which are
only required to 
be \emph{nonanticipative} and allow service 
\emph{preemptions} at the start of
each  period.

The system incurs linear holding and/or
 rejection costs separably across classes, time-discounted with factor
 $0 < \beta < 1$. 
Class $k$ traffic incurs holding costs at rate $c_k \geq 0$ per period
 and customer in the system, and rejection
 costs at rate $r_k \geq 0$ per customer blocked.
Denoting by $A_k \sim \text{Poisson}(\lambda_k)$
  the number of class
 $k$ arrivals during a period,  we can thus represent the 
class' state-dependent cost rate function as follows: 
for $(a_k, i_k) \in N_k$,
\begin{equation}
\label{eq:hkjk}
h_k(a_k, i_k) = h_k(i_k) \triangleq 
c_k i_k + r_k \mathbb{E}\left[\left\{A_k - (n_k - i_k)\right\}^+\right].
\end{equation}

The goal is to design a tractable index policy that comes close to
minimizing the expected total discounted value of costs incurred.
Since the model is readily formulated as a restless bandit problem,
the natural candidate for such an index is the MPI.
Note that the latter is a function $\nu_k^*(1, i_k)$ of the controllable
states of each class $k$. 

Along the lines in \citet{nmqs06} (cf.\ Sect.\ \ref{s:smqfb})
we distinguish between two types of traffic classes.
We consider a class $k$ to be \emph{loss-sensitive} if 
$(1-\beta) r_k \geq c_k - \epsilon_k$, and \emph{delay-sensitive} if 
$(1-\beta) r_k \leq c_k - \epsilon_k$, where $\epsilon_k > 0$ 
is a critical value whose characterization is discussed in
 the full version of \citet{nmngi06}, currently under preparation.

In short, the indexability analyses reveals existence of the MPI
and structural properties that extend those in the fully-connected
special case in \citet{nmqs06}.
Thus, the ordering induced by the MPI $\nu_k^*(1, i)$ on the state
space of a class depends on whether this is loss-sensitive or
delay-sensitive, with such an ordering being in each case the same as 
that discussed in Sect.\ \ref{s:smqfb} for the corresponding MPI
$\nu_k^*(i)$.

Detailed analyses will be given in the full version of
\citet{nmngi06}, along with comprehensive experimental results
assessing the performance of such an MPI policy.

\subsection{Scheduling a multiclass queue with finite buffers and
     delayed state observation}
\label{s:smqfbdso}
In some applications, it appears unreasonable due to communication
delays to assume that the
controller has knowledge of the
current system state, which has motivated research into control of
systems with delayed state information. See, e.g., \citet{altsti} and the
references therein.
In the setting of communication systems for scheduling multiclass
traffic, such delays are significant in, e.g., satellite systems. 
This has motivated recent research efforts to address the design of 
dynamic scheduling policies in multiclass queues with delayed state
observation.
Thus, \citet{ehsanliu06} establish the optimality of a
simple greedy scheduling policy, yet under strong conditions on model parameters.

In \citet{nmngi07} we announce results of deploying restless bandit
indexation in a multiclass finite-buffer queueing scheduling model which
differs from that discussed in Sect.\ \ref{s:smwqfb} in that: (i)
classes are permanently connected to the server; and (ii) the
controller's state information on queues' states is delayed by one
time period.
The model is readily cast as a restless bandit problem, where the
\emph{observed state} process of each bandit (class) $k$ is of the
form $\widehat{X}_k(t) \triangleq (a_k(t-1), X_k(t-1))$, i.e., it
consists of the previous action $a_k(t-1) \in \{0, 1\}$ and buffer
occupancy $X_k(t-1)$.
Note that in that setting all states for a class are controllable,
as the controller might well allocate the server to a class whose
previous backlog was empty, anticipating that such need not be the case in the
present period.

Again, we find that the required family $\widehat{\mathcal{F}}$ of
active sets relative to which indexability is established depends
critically on whether a class is loss-sensitive or delay-sensitive. 
The former case corresponds to a class $k$ satisfying, in the notation
of Sect.\ \ref{s:smwqfb}, that $0 < (1-\beta) r_k \geq \beta c_k$, whereas
the latter case corresponds to satisfaction of $(1-\beta) r_k \leq
\beta c_k > 0$.
The paper discusses the structure of $\widehat{\mathcal{F}}$ and of
the MPI that emerges in each case, and further reports results of
preliminary computational experiments, where the MPI policy is shown
to be nearly optimal and to outperform conventional scheduling
policies.
A full version of the paper with detailed analyses and thorough
computational experiments is currently under preparation.

\subsection{Multiarmed bandits with switching costs}
\label{s:mabsc}
A critical assumption underlying the optimality of the Gittins index
rule for the multiarmed bandit problem is that switching projects is
costless.
Yet, in many applications there are nonnegligible costs of switching,
 motivating investigation of multiarmed bandits with switching
costs, which are extensively surveyed in \citet{jun04}.
In such a setting, it is natural to consider index policies where the 
index $\nu_{(a^-, i)}$ of a project is a function of its present state
$i$ and previous action $a^-$ (i.e., $a^- = 1$ if it was active and $a^- =
0$ otherwise).
As pointed out in  \citet{banksun94}, the presence of switching costs
leads one to stick longer to the project currently engaged. Such a 
\emph{hysteresis} property, also discussed in \citet{duhon03} in the
setting of a continuous-state model, means that the indices should be
required to satisfy $\nu_{(1, i)} \geq \nu_{(0, i)}$, i.e., the 
\emph{continuation index} is larger than or equal to the 
\emph{switching index}.

Though \citet{banksun94} established that, generally, optimal policies
for bandits with switching costs are not of index type, 
\citet{asatene96} introduced an intuitively appealing index, and went
on to prove that it partially characterizes optimal policies.
Their  continuation index $\nu_{(1, i)}^{\textup{AT}}$ is precisely
the project's Gittins index, whereas their switching index
$\nu_{(0, i)}^{\textup{AT}}$ is the maximum rate of expected
discounted reward minus initial switching cost per unit of expected
discounted time that can be achieved under stopping rules that engage
an initially passive project.

\citet{nmswc07} deploys restless bandit indexation to such a problem,
by exploiting the natural restless reformulation of a classic bandit with
switching costs, where one considers the augmented state $(a^-,
i)$.
As the appropriate family $\widehat{\mathcal{F}}$ of active sets in
the augmented state space $\widehat{N} = \{0, 1\} \times N$, where $N$
is the project's original state space, we take 
\begin{equation}
\label{eq:fhat}
\widehat{\mathcal{F}} \triangleq \left\{S_0 \oplus S_1\colon S_0
  \subseteq S_1 \subseteq N\right\},
\end{equation}
where the notation $S_0 \oplus S_1$ refers to the policy that engages
the project when it was previously rested (resp.\ engaged) iff its
current state lies in $S_0$ (resp.\ $S_1$).

The paper shows that the resultant restless bandits are
PCL$(\widehat{\mathcal{F}})$-indexable, and that their MPI is
precisely the Asawa and
Teneketzis index. 
More importantly, computational issues are addressed. 
\citet{asatene96} had proposed to jointly compute both the active and
the passive index of an $n$-state project as the Gittins index of a
certain classic project having $2 n$ states. This results in an
eight-fold increase in $O(n^3)$ arithmetic operations and
substantially increased memory operations relative to the
computational effort to compute the continuation index alone, which
is overly expensive in large-scale models.

We thus set out to obtain a more efficient computation method by
analyzing the adaptive-greedy algorithm $\mathrm{AG}_{\widehat{\mathcal{F}}}$
in such a model. It turns out that the algorithm naturally decouples
into a two-stage method. The first stage computes the Gittins index of
the original project along with certain required extra quantities. 
Then, the second stage is fed the first stage's output to compute the 
switching index \emph{an order of magnitude faster} in at most 
$n^2 + O(n)$ arithmetic operations. The two-stage method also yields 
substantially reduced memory operations, as it involves manipulation
of $n \times n$ matrices instead of $2 n \times 2n$ matrices. 

A computational study demonstrates that such a theoretical complexity
reduction translates in practice into dramatic runtime savings. 
Further, the paper reports on a comprehensive computational study on
two- and three-class project instances showing that the MPI policy is
consistently near optimal, and substantially outperforms the Gittins index
policy that ignores switching costs.

To help the reader grasp the underlying intuition, 
 we present next
an illustrative example.  
Consider the 3-state bandit instance with constant
startup cost $c$,
\[
\beta = 0.95, \quad  \mathbf{R} = \begin{bmatrix} 0.0250 \\ 0.4242 \\ 0.0338\end{bmatrix}, 
\quad \text{ and } \quad 
\mathbf{P} = 
\begin{bmatrix} 
0.6635 & 0.0285 & 0.3080 \\
0.6345 & 0.3583 & 0.0072 \\
0.4868 & 0.0530 & 0.4602
\end{bmatrix}.
\]

\begin{figure}[htb]
\centering
\begin{psfrags}%
\psfragscanon%
%
\psfrag{s03}[t][t]{\color[rgb]{0,0,0}\setlength{\tabcolsep}{0pt}\begin{tabular}{c}$g^{\pi}$\end{tabular}}%
\psfrag{s04}[b][b]{\color[rgb]{0,0,0}\setlength{\tabcolsep}{0pt}\begin{tabular}{c}$f^{\pi}$\end{tabular}}%
\psfrag{s05}[b][b]{\color[rgb]{0,0,0}\setlength{\tabcolsep}{0pt}\begin{tabular}{c}no startup cost\end{tabular}}%
\psfrag{s06}[t][t]{\color[rgb]{0,0,0}\setlength{\tabcolsep}{0pt}\begin{tabular}{c}$g^{\pi}$\end{tabular}}%
\psfrag{s07}[b][b]{\color[rgb]{0,0,0}\setlength{\tabcolsep}{0pt}\begin{tabular}{c}$f^{\pi}$\end{tabular}}%
\psfrag{s08}[b][b]{\color[rgb]{0,0,0}\setlength{\tabcolsep}{0pt}\begin{tabular}{c}startup cost $c = 0.02$\end{tabular}}%
%
\includegraphics[height=2.2in,width=6.5in,keepaspectratio]{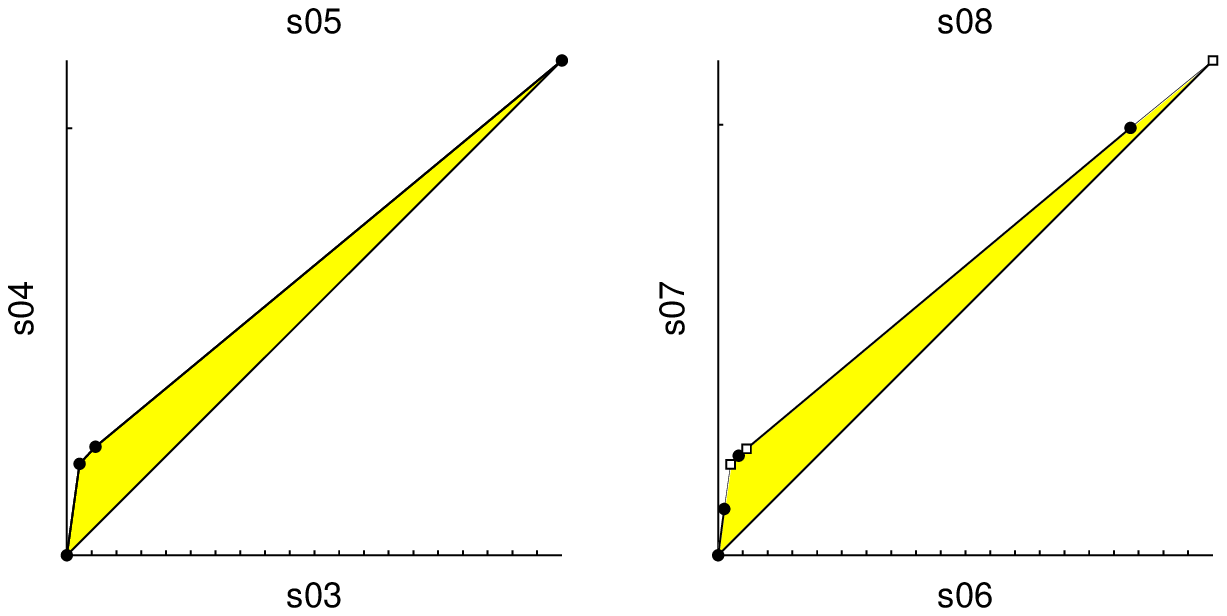}
\end{psfrags}%
\caption{Achievable work-reward performance regions and structure of upper boundaries.}
\label{fig:ar1}
\end{figure}

Work and reward measures $g^\pi$ and $f^\pi$ are evaluated 
assuming that the initial state is uniformly drawn.
The left pane in Figure \ref{fig:ar1} shows the 
achievable work-reward performance region in the 
classic $c = 0$ case.
The four points displayed, determining its \emph{upper boundary}, are the
work-reward performance points corresponding to the policies whose
active sets are given, from left to right, by
$\emptyset$, $\{2\}$, $\{2, 3\}$, and 
$\{2, 3, 1\}$.
The successive work-reward trade-off slopes or rates between such points are
the Gittins index values:
\[
\nu_2^* = 0.4242 > \nu_3^* = 0.061487 > 
\nu_1^* = 0.048002.
\]

The right pane in Figure \ref{fig:ar1} shows a corresponding plot for
the case with startup cost $c = 0.02$.
The upper work-reward boundary is determined by 
the seven points displayed, which are the 
work-reward performance points corresponding, from left to right, to
the policies having active sets
$\emptyset \oplus \emptyset$,
$\emptyset \oplus \{2\}$,
$\{2\} \oplus \{2\}$, 
$\{2\} \oplus \{2,  3\}$,
$\{2, 3\} \oplus \{2, 3\}$, $\{2, 3\} \oplus \{2, 3,
1\}$
 and  $\{2, 3, 1\} \oplus \{2, 3,
1\}$.
The successive work-reward trade-off slopes between such points give
the MPI values:
\[
\nu_{(1, 2)}^* = 0.424 > \nu_{(0, 2)}^* = 0.411 > \nu_{(1, 3)}^* = 0.061 > 
\nu_{(0, 3)}^* = 0.051 > \nu_{(1, 1)}^* = 0.048 > \nu_{(0, 1)}^* = 0.047.
\] 
The plot represents the right end-points giving a
continuation index value by a black circle, and those giving a
switching index value by a white square.
Notice further that the continuation index matches the Gittins
index of the previous case.

\begin{figure}
\centering
\begin{psfrags}%
\psfragscanon%
\psfrag{s01}[t][t]{\color[rgb]{0,0,0}\setlength{\tabcolsep}{0pt}\begin{tabular}{c}$g^{\pi}$\end{tabular}}%
\psfrag{s02}[b][b]{\color[rgb]{0,0,0}\setlength{\tabcolsep}{0pt}\begin{tabular}{c}$f^{\pi}$\end{tabular}}%
\psfrag{s03}[b][b]{\color[rgb]{0,0,0}\setlength{\tabcolsep}{0pt}\begin{tabular}{c}startup cost $c = 0.1$\end{tabular}}%
\psfrag{s05}[t][t]{\color[rgb]{0,0,0}\setlength{\tabcolsep}{0pt}\begin{tabular}{c}$g^{\pi}$\end{tabular}}%
\psfrag{s06}[b][b]{\color[rgb]{0,0,0}\setlength{\tabcolsep}{0pt}\begin{tabular}{c}$f^{\pi}$\end{tabular}}%
\psfrag{s07}[b][b]{\color[rgb]{0,0,0}\setlength{\tabcolsep}{0pt}\begin{tabular}{c}startup cost $c = 0.6$\end{tabular}}%
%
\psfrag{x01}[t][t]{0}%
\psfrag{x02}[t][t]{0.1}%
\psfrag{x03}[t][t]{0.2}%
\psfrag{x04}[t][t]{0.3}%
\psfrag{x05}[t][t]{0.4}%
\psfrag{x06}[t][t]{0.5}%
\psfrag{x07}[t][t]{0.6}%
\psfrag{x08}[t][t]{0.7}%
\psfrag{x09}[t][t]{0.8}%
\psfrag{x10}[t][t]{0.9}%
\psfrag{x11}[t][t]{1}%
%
\psfrag{v01}[r][r]{0}%
\psfrag{v02}[r][r]{0.1}%
\psfrag{v03}[r][r]{0.2}%
\psfrag{v04}[r][r]{0.3}%
\psfrag{v05}[r][r]{0.4}%
\psfrag{v06}[r][r]{0.5}%
\psfrag{v07}[r][r]{0.6}%
\psfrag{v08}[r][r]{0.7}%
\psfrag{v09}[r][r]{0.8}%
\psfrag{v10}[r][r]{0.9}%
\psfrag{v11}[r][r]{1}%
%
\includegraphics[height=2.2in,width=6.5in,keepaspectratio]{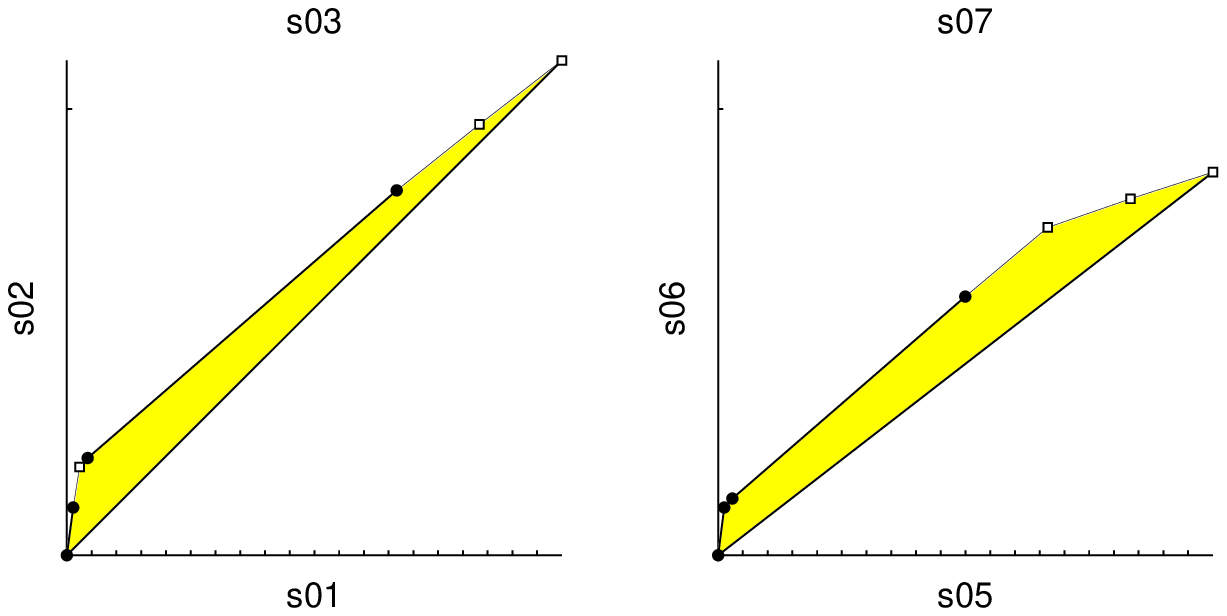}
\end{psfrags}%
\caption{Achievable work-reward performance regions and structure of upper boundaries.}
\label{fig:ar2}
\end{figure}

The left pane of Figure \ref{fig:ar2} shows the achievable
work-reward performance region for the case $c = 0.1$. 
Now, the seven points displayed in the upper boundary correspond, from left to right, to the
policies having active sets
$\emptyset \oplus \emptyset$, 
$\emptyset \oplus \{2\}$,
$\{2\} \oplus \{2\}$, 
$\{2\} \oplus \{2, 3\}$, 
$\{2\} \oplus \{2, 3, 1\}$,
$\{2, 3\} \oplus \{2, 3, 1\}$
and 
$\{2, 3, 1\} \oplus \{2, 3, 1\}$.
The successive work-reward trade-off slopes give
the bandit's MPI values:
\[
\nu_{(1, 2)}^* = 0.424 > \nu_{(0, 2)}^* = 0.358 > \nu_{(1, 3)}^* =
0.061 > \nu_{(1, 1)}^* =
0.048 > \nu_{(0, 3)}^* = 0.044 > \nu_{(0, 1)}^* = 0.043.
\] 

Finally, the right pane of Figure \ref{fig:ar2} shows the
corresponding plot for the case $c = 0.6$.
The seven points characterizing the upper boundary
correspond, from left to right, 
to the
policies having active sets $\emptyset \oplus \emptyset$, 
$\emptyset \oplus \{2\}$,
$\emptyset \oplus \{2, 3\}$,
$\emptyset \oplus \{2, 3, 1\}$,
$\{2\} \oplus \{2, 3, 1\}$,
$\{2, 3\} \oplus \{2, 3, 1\}$, and
$\{2, 3, 1\} \oplus \{2, 3, 1\}$.
The resultant MPI values given by the successive slopes are
\[
\nu_{(1, 2)}^* = 0.424 > \nu_{(1, 3)}^* = 0.061 > \nu_{(1, 1)}^* =
0.048 >  \nu_{(0, 2)}^* = 0.047 > \nu_{(0, 3)}^* = 0.019 > \nu_{(0, 1)}^* = 0.018.
\] 

Notice that, in each case, the continuation index value $\nu_{(1,
  i)}^*$ matches the Gittins index value $\nu_i^*$.
Further, the successive active sets $S_0 \oplus S_1$ 
characterizing the regions' upper boundaries belong in the
active-set family $\widehat{\mathcal{F}}$ in (\ref{eq:fhat}).
Also, the continuation index value $\nu_{(1, i)}^*$ is 
larger than the corresponding switching index value $\nu_{(0,
  i)}^*$ value.

\subsection{Multiarmed bandits with switching delays}
\label{s:mabsd}
Besides switching costs, delays for switching projects are clearly
relevant in many applications. Think, e.g., of the time required to
learn a new technique before one can make productive use of it, of
time for preparing the ground in a development project, or of the
shutdown time required to dismantle a project.
\citet{asatene96} briefly discussed the corresponding multiarmed bandit
problem  with switching delays, for which they proposed an index,
claiming that it also gives a partial characterization of optimal
policies. Yet, no algorithm is given in their paper for the
computation of such an index.

\citet{nmswd07} announces results of extending the restless bandit indexation analysis in
\citet{nmswc07} to projects that incorporate both switching costs and
delays.
The required extension is, however, far from trivial. First, the
natural restless reformulation of a classic project with switching
delays is semi-Markovian. Further and more importantly, the resultant
restless projects are no longer PCL-indexable. 
In fact, it was the analysis of this model which motivated the author
to develop the more powerful LP-indexability conditions reviewed in
Sect.\ \ref{s:mpicfic}.

Thus, the paper shows that the restless projects of concern are
LP$(\widehat{\mathcal{F}})$-indexable, where $\widehat{\mathcal{F}}$ is the active-set
family discussed in the previous section. 
This allows application of the adaptive-greedy algorithm
$\mathrm{AG}_{\widehat{\mathcal{F}}}$ to compute the MPI. 
Similarly as in the case of switching costs only, analysis of such an
algorithm as it applies to the model with switching delays reveals
that it can be decoupled into a two-stage scheme. 
Again, the first stage computes the Gittins index of the original
$n$-state project along with extra quantities, and the second stage
uses the first stage's output to compute the switching index. 
This is accomplished an order of magnitude faster, in at most $(5/2)
n^2 + O(n)$ arithmetic operations.

The paper further reports on a computational study showing
the dramatic speedup gains yielded by such a two-stage computation
scheme, relative to joint computation of both indices.
Such a study is complemented by a set of experiments aimed at
assessing the degree of suboptimality of the MPI policy, and its
performance gains over the benchmark index policy that ignores
switching penalties. These experiments reveal substantial gains and a
near-optimal performance.

\subsection{Multiarmed bandits with deadlines}
\label{s:mabd}
Another critical model assumption underlying the Gittins and Jones
result for classic multiarmed bandits is that the planning horizon is
infinite. 
Yet, in many applications it is more appropriate to consider finite
horizon scenarios. 
Thus, projects might be subject to a common deadline at which they expire;
or, more generally, each project might have
its own deadline, which motivates consideration of the 
\emph{multiarmed bandit problem with deadlines}.

While such problems appear to be
computationally intractable, their widespread practical relevance motivates the
quest for tractable priority-index policies that are nearly optimal.
For such a purpose, the index introduced in \citet{bjk56}, as it
extends to a general Markovian setting, is a natural choice. 
Such an index measures the maximum rate of expected discounted reward
per unit of expected discounted time that can be achieved under
stopping rules that do not exceed the horizon.
\citet{gijo74} discussed the resultant index rule, showing that
it is generally not optimal for the finite-horizon multiarmed bandit
problem.

Remarkably, to the best of the author's knowledge, use of such a priority-index rule
has received scant, if any, research attention.
It appears that the finite-horizon index has been regarded in a
subordinate role, simply as a means to approximate the Gittins index. 
See, e.g., \citet[Sect.\ 7]{gi79}, \citet[Chapt.\ 1]{gi89}, and
\citet{wang97}.

One possible explanation accounting for such a limited use of the
finite-horizon index might be the lack of a simple, exact algorithm
for its computation. 
The author is only aware of the exact method for 
Bayesian Bernoulli bandits with beta priors discussed
in  \citet[Sect.\ 7]{gi79}. Yet,
such a method involves the numerical computation of a supremum at
each step --- see formula (11) in that paper --- which is not an elementary
operation.

An approach to address such issues based on restless bandit indexation
is announced in \citet{nmcdc05}, based on the observation that  a
finite-horizon bandit is readily formulated as an infinite-horizon
restless bandit, by augmenting the state to include the remaining
time, and on the result that the latter's MPI is the
former's finite-horizon index.
In this setting, it is natural to represent Markov deterministic
policies for operating a project in the form
\[
\widehat{S} = S_1 \oplus S_2 \oplus \cdots S_T,
\]
where the notation represents the policy that engages the project
when $t$ periods remain to the deadline iff the original state lies in
$S_t$, for $1 \leq t \leq T$.
The restless bandits of concern for a $T$-horizon
project turn out to be
PCL$(\widehat{\mathcal{F}}_T)$-indexable, where 
 the appropriate active-set family $\widehat{\mathcal{F}}_T$
is given by 
\[
\widehat{\mathcal{F}}_T \triangleq \left\{S_1 \oplus S_2 \oplus \cdots
  \oplus S_T\colon S_1 \subseteq S_2 \subseteq \cdots \subseteq S_T\right\}.
\]
Notice that such a structure is consistent with the intuitive
monotonicity property of the index whereby $\nu_{(t, i)}^* \leq
\nu_{(t-1, i)}^*$, with $t$ being the remaining time and $i$ the
current state.

Such a result allows us to use the adaptive-greedy algorithm
$\mathrm{AG}_{\widehat{\mathcal{F}}_T}$ to compute the MPI. However,
the resultant complexity is of order $O(T^3 n^3)$ arithmetic
operations for a $T$-horizon $n$-state project, which severely hinders applicability.
In more recent work, we improve on such a result, by 
showing that the adaptive-greedy algorithm naturally decouples into a
recursive $T$-stage scheme, which computes the finite-horizon MPI 
in only $O(T^2 n^3)$ arithmetic operations, thus significantly expanding the size of models
that can be addressed.
Finally, preliminary computational  results on
two-project instances show that such MPI policies are consistently
near optimal, often substantially outperforming the benchmark Gittins
index policy.

\section{Concluding remarks}
\label{s:concl}
We have surveyed a unifying approach to design and compute tractable
priority-index policies for a variety of problems, based on the
intuitive concept of MPI.
We believe that, while the results reviewed on theory and algorithms
for finite-state bandits are in a state that allows them to be readily
deployed by researchers, the applications surveyed only give a glimpse
of what can be attained both in depth and scope.
Many interesting issues remain to be explored, such as: Is the
prevalence of indexability and PCL-indexability as high as it seems in
computational experiments,
and, if so, Why?
How to analyze the performance of MPI policies in a
multi-project setting? 
Under what conditions do MPI policies perform well?
How can one design suitable MPI policies in models with constraints,
as in \citet{altscw89}?
We emphasize that the range of applications of the MPI
approach is far from exhausted. In fact, the author is currently
investigating completely different applications from those surveyed
herein, which will be reported in due time.
Also, Prof.\ Weber has introduced in his discussion to this paper a promising idea
that might substantially expand the scope of the MPI approach.

\subsection*{Acknowledgements}
The author's work surveyed in this paper
has been supported in part by the Spanish Ministry
of Education \& Science under projects TAP98-0229, BEC2000-1027, 
 MTM2004-02334, a Ram\'on y Cajal Investigator Award and an I3 faculty
 endowment grant,
by the European Union's Networks of Excellence 
  EuroNGI and EuroFGI, 
and by the Autonomous Community of Madrid under grants
  UC3M-MTM-05-075 and
CCG06-UC3M/ESP-0767.
The author thanks, for their invitations to present parts of his work
surveyed herein, the organizers of research seminars at 
Eurandom (1999, 2001, 2004) and at Beta (2004), Technical Univ.\ of Eindhoven, and the
Department of Operations Research and Statistics at Univ.\ de
Sevilla (2007), as well as the organizers of the 
Workshop on Fluid Queues (Eurandom, 2000), 
the Schloss Dagstuhl Seminar on
 Scheduling in Computer and Manufacturing Systems (Wadern, Germany,
 2002), 
the Workshop on Analysis and Optimization of
Stochastic Networks with Applications to Communications and
Manufacturing (Eurandom, 2002), 
the Workshop on 
Topics in Computer Communication and Networks (Isaac Newton Institute for Mathematical Sciences, 
Univ.\ of Cambridge, UK, 2002),
the Joint Research Conference on Mathematical Programming 
of the Spanish Operations Research
   \& Statistics Society
 and the Spanish Royal Mathematical Society (Univ.\ Miguel
 Hern\'andez, Elche,
Spain, 2004), 
the Schloss Dagstuhl Seminar on
 Algorithms for Optimization with Incomplete Information (Wadern,
 Germany, 2005), and the 
First Iberian Conference on Optimization (Coimbra, Portugal, 2006).

\bibliographystyle{spbasic}      

\end{document}